%% file: ST_KM_0110.tex
\documentclass{svjour3}                     % onecolumn (standard format)
\smartqed  % flush right qed marks, e.g. at end of proof
\usepackage{graphicx}
\usepackage{amsmath}
\usepackage[T1]{fontenc}
\usepackage{amsfonts}
\usepackage[english]{babel}
\usepackage[dvips]{color}
\usepackage{amsmath,amssymb}

\def\func#1{\rm {#1}}

\def\cE{{\mathcal{E}}}

\def\mb0{{\mathbf{0}}}

\def\leftB{[\![}
\def\rightB{]\!]}
\include{macro}

\begin{document}

\title{Weak error for stable driven SDEs  
: expansion of the densities.%\thanks{Grants or other notes
}
%\subtitle{Euler scheme for stable SDEs}

\titlerunning{Euler scheme for stable SDEs}        % if too long for running head

\author{Valentin Konakov \and St\'ephane Menozzi %etc.
}

\authorrunning{V. Konakov and S. Menozzi} % if too long for running head

\institute{V. Konakov, \at
              CEMI, Academy of Sciences. Nahimovskii av. 47, 117418 Moscow, RUSSIA. \\
              \email{valentin\_konakov@yahoo.com}           %  \\
           \and
           S. Menozzi, \at
              LPMA, Universit\'e Paris VII Diderot.  175 Rue du Chevaleret, 75013 Paris, FRANCE.\\
              \email{menozzi@math.jussieu.fr}  
}

%\institute{V. Konakov 
%\at CEMI RAS, Moscow
%%Central economics Mathematical institute,\\
%%Academy of sciences,\\
%%Nahimovskii av. 47,\\
%%117418 MOSCOW.\\
%%RUSSIA.\\        
%              \email{valentin\_konakov@yahoo.com}          \\
%%             \emph{Present address:} of F. Author  %  if needed
%           S. Menozzi \at
%%Laboratoire de Probabilit\'es\\
%%et Mod\`eles Al\'eatoires\\
%%Universit\'e Paris VII,\\
%%175 Rue du Chevaleret, 75013 PARIS.\\
%%FRANCE.\\
%               
%                \email{menozzi@math.jussieu.fr}  
%}

\date{\today}

\maketitle

\begin{abstract}
Consider a multidimensional SDE of the form
$X_t=x+\int_{0}^{t}b(X_{s-})ds+\int_{0}^{t}f(X_{s-})dZ_s $ where
$(Z_s)_{s\ge 0}$ is a symmetric stable process. %(possibly perturbed by a compound Poisson process). 
Under suitable
assumptions on the coefficients the unique strong solution of the
above equation admits a density w.r.t. the Lebesgue measure and so
does its Euler scheme. Using a parametrix approach, we derive an
error expansion %at order 1 
w.r.t. the time step for the difference
of these densities.
\keywords{Symmetric stable processes\and  parametrix \and Euler scheme  }
%\PACS{PACS code1 \and PACS code2 \and more}
\subclass{60H30\and 65C30 \and 60G52 }
\end{abstract}

 \mysection{Introduction}
Consider the following $\R^d $-valued 
Stochastic Differential Equation (SDE in short)
\begin{equation} 
X_t=x+\bint{0}^{t}b(X_{s-})ds+\bint{0}^t f(X_{s-}) dZ_s,
 \label{GEN_EDS}
\end{equation}
where $b, f $ are respectively Lipschitz continuous mappings from $\R^d $ to $\R^d $ and $\R^d $ to $\R^d\otimes  \R^d$ and $(Z_s)_{s\ge 0}$ is a general L\'evy process. The previous assumptions guarantee the existence of a unique strong solution to  \eqref{GEN_EDS}. Also, this solution satisfies the strong Markov property, see e.g. Theorem 7 and 32 Chapter 5 in Protter \cite{prot:04}. 
Let $T>0$ be a fixed time horizon and  $(X_{t}^N)_{t\in\Lambda}$ a given approximation scheme of $(X_t)_{t\in[0,T]} $ associated to the time step $h:=T/N,\ N\in\N^* $ on the grid $\Lambda:=\{t_i:=ih,i\in\leftB 0,N\rightB \}$. When speaking about weak approximation of \eqref{GEN_EDS} two kinds of quantities are of interest. The first one writes 
\begin{eqnarray*}
{\cal E}_1(x,T,N):=\E_x[g(X_T)]-\E_{x}[g(X_T^N)]
\end{eqnarray*}
for a suitable class of test functions $g$.  The second one concerns, when it exists, the approximation of the  transition density $p $ of the original SDE \eqref{GEN_EDS}. If the approximation scheme $(X_t^N)_{t\in \Lambda} $ admits as well a transition density $p^N$, the quantity under study becomes 
\begin{eqnarray*}
{\cal E}_2(x,y,T,N):=(p-p^N)(T,x,y).
\end{eqnarray*}
In both cases, the goal is to give a bound or an error expansion of these quantities in terms of $h$. The error expansions are particularly useful for practical simulation. For ${\cal E}_1$, the expansion allows to use the Romberg Richardson extrapolation to improve the convergence of the discretization error, see e.g. Talay and Tubaro \cite{tala:tuba:90}. 
On the other hand, if $p $ and $p^N$ exist, and a suitable expansion of ${\cal E}_2  $ holds, %, one can use it to get a control on ${\cal E}_1 $ for a wider class of tests functions $g$. Namely, the growth assumptions on $g$ can be closely related to the decay of the density at infinity. Also, a precise development of ${\cal E}_2$ can be useful 
it can be useful to estimate the sensitivity of ${\cal E}_1 $ w.r.t. to the spatial variable $x$ and it also allows to get a control on ${\cal E}_1 $ for a wider class of test functions $g$ than those considered by the direct methods used to control this quantity, see e.g. Guyon \cite{guyo:06}.  
Indeed, the typical assumptions and techniques associated to the 
study of $ \cE_1$ and $\cE_2 $ are quite of different nature.

In the continuous case, i.e. $Z_s=bs+\sigma W_s$ where $(W_s)_{s\ge 0} $ is a standard $d$-dimensional Brownian motion, provided the test function $g$ and the coefficients $b,f $ are sufficiently smooth and $g$ has polynomial growth, without any additional assumption on the generator Talay and Tubaro \cite{tala:tuba:90} derive an error expansion at order 1 for $\cE_1(x,T,N)$ when $(X_t^N)_{t\in\Lambda}$ is the Euler approximation. Their proof is based on standard stochastic analysis tools: It\^o's expansions and stochastic flows. To obtain the same kind of result for bounded Borel functions $g$ some non degeneracy has to be assumed, namely hypoellipticity of the underlying diffusion, and the proof relies on Malliavin calculus techniques, see Bally and Talay \cite{ball:tala:96:1}. The authors also manage to extend their results to $\cE_2(x,y,T,N)$ for a slightly modified Euler scheme \cite{ball:tala:96:2}.

Anyhow, in the uniformly elliptic case, the most natural approach to handle the estimation of  the quantity $\cE_2(x,y,T,N) $ consists in using the so called "parametrix" technique introduced to obtain existence and controls on the fundamental solutions of PDEs, see e.g. Mc Kean and Singer \cite{mcke:sing:67} or Friedman \cite{frie:64}. Roughly speaking it consists in expressing the density of $ X_T$ in terms of an infinite sum of suitable iterated kernels applied to the density of an SDE with constant coefficients.
This has been done successfully by Konakov and Mammen \cite{kona:mamm:02}.% in a uniformly elliptic setting. 
The main advantage of this approach is that the density of the solution $X_T$ and the Euler approximation $X_T^N$ can be expressed in the same form and therefore quite directly compared. Furthermore this technique turns out to be quite robust and can be applied as soon as good controls on the densities $p,p^N$ and their derivatives are available, see e.g. \cite{kona:meno:molc:08} for an extension to a slightly degenerate framework. 

%Consider now the general case, i.e. for all $u\in\R^d $, $\E[\exp(i\langle u,Z_t\rangle)]=\exp\left(t\psi(u)\right) $ with
%$$ \psi(u)=\left\{i\langle \gamma ,u\rangle-\sigma^2\frac{|u|^2}{2}+\bint{}^{}\left(\exp(i\langle u,z\rangle)-1-i\langle u,z\rangle\I_{|z|\le 1} \right)\nu(dz)\right\},$$
%where $\gamma\in\R^d,\ \sigma  \in \R$ and $\nu$ stands for the L\'evy measure of the process.
For  a general L\'evy process $Z$ and suitable smooth functions $b,f$, $g$, under additional assumptions on the behavior at infinity of the L\'evy measure $\nu $ of $Z$, that is integrability conditions of the large jumps, Protter and Talay \cite{prott:tala:97}, manage to get a control at order one or even an error expansion for $\cE_1(x,T,N)$  with the same approach as in \cite{tala:tuba:90}. In that work the approximation is the Euler scheme which for a general L\'evy measure $\nu $ cannot always be exactly simulated on a computer. %%% Regarder le resultat de Hausenblas

The quantity $\cE_1(x,T,N) $ for approximations of the Euler scheme that can be simulated
has also been studied by Jacod {\textit {et al.}} \cite{jaco:etal:05} who derived bounds at order 1. Moment conditions are also assumed. We finally refer to the work of Hausenblas and Marchis \cite{haus:marc:06} for approximations of Poisson jump measures that are easy to simulate.
\\

In this work, we consider the case where $(Z_t)_{t\ge 0}$ is an $\alpha $-stable symmetric process, $\alpha\in(0,2)$.
Under suitable non degeneracy assumptions on its coefficients specified below (see \A{A-1}-\A{A-3}),
equation \eqref{GEN_EDS} is known to have a density $p$ w.r.t. the
Lebesgue measure. This can be proved via a Malliavin calculus-Bismut integration by parts
approach, see e.g. Bichteller \textit{et al.}
\cite{bich:grav:jaco:87}. Also, 
a direct construction of this
density using a parametrix expansion has been obtained by
Kolokoltsov \cite{kolo:97} who derived as well "Aronson's like"
bounds with time singularity depending on the index $\alpha $ of the
stable process $(Z_t)_{t\ge 0}$. % and the tails of the compound Poisson process $Z^P$.

Analogoulsy to the "diffusion case" the first step of the parametrix is to consider that the density $p(T,x,y) $ of \eqref{GEN_EDS} can be approximated by the density of the process $\widetilde X_t^y=x+b(y)t+f(y)Z_t$ at time $T$. Namely, we freeze the coefficients in \eqref{GEN_EDS} at the final spatial point. The next crucial point is to obtain sharp estimates  of the stable density $\widetilde p^y(T,x,.) $ of $\widetilde X_T^y $ and its derivatives in order to solve the parametrix integral equations.

%%%%%%%%%%%% Ce que l'on metnionne au paragraphe precedent
% Pour resumer les deux grands types d'hypotheses en jeu.
% BILAN REGULARITE ET MOMENTS, OU NON DEGENERESECENCE en un certain sens 

Stable driven SDEs appear in various applicative fields, from
mathematical physics to electrical engineering or financial mathematics, see \cite{imke:pavl:06}, \cite{stuc:klei:74} or \cite{jani:mich:wero:96}, therefore their
approximation becomes of interest. To
approximate equation \eqref{GEN_EDS}, setting 
$\phi(t):=\inf\{t_i: t_i\le t<t_{i+1} \}$, we introduce the Euler scheme
\begin{equation}
\label{EUL_SCH} X_t^N=x+\bint{0}^{t}b(X_{\phi(s)}^N)ds
+\bint{0}^{t}f(X_{\phi(s)}^N)dZ_s.
\end{equation}
The computation of the above scheme only requires to be able to simulate exactly the increments  
of $(Z_t)_{t\ge 0}$, which up to a self similarity argument only amounts to simulate a stable law. This aspect is for instance discussed in Samorodnitsky and Taqqu \cite{sam:taqq:94}, Weron and Weron \cite{wero:wero:95} or Section 3 of \cite{prott:tala:97}.
Under the same assumptions \A{A-1}-\A{A-3}, the Euler scheme
defined in \eqref{EUL_SCH} also has a density $p^N$.

Observe that the results of \cite{prott:tala:97}, \cite{jaco:etal:05} cannot be directly applied, even for the study of $\cE_1(x,T,N) $, since stable laws have heavy tails. Comparing the parametrix developments of $p$ and $p^N$ we obtain an expansion with leading
term of order 1 in $h$ for $\cE_2(x,y,T,N) $. The
parametrix expansion of $p$ is discussed in \cite{kolo:97}, see also Section \ref{CONN_SDE_IDO} and Appendix, whereas the
parametrix expansion of $p^N$ can be related to the ideas developed
in \cite{kona:mamm:00,kona:mamm:02} for the diffusive case
corresponding to an index of stability equal to $2$.

This result also emphasizes the robustness of the method that naturally extends to a broad class of processes. Let us mention that, using a Malliavin calculus approach, Hausenblas \cite{haus:02}, derived an upper bound of order one w.r.t. $h$ for the quantity $\cE_1(x,T,N),\ g\in L^\infty$ in the scalar case. Concerning functional limit theorems for the approximation of stable driven SDEs we refer to the work of Jacod \cite{jaco:04}. 

The paper is organized as follows. In Section \ref{ASS_and_RES} we state our standing assumptions and main results. In Section \ref{CONN_SDE_IDO} we prove the existence of the densities for both the  stable driven equation and its Euler scheme and also give a parametrix representation of these densities.
Section \ref{PROOFS_GB} is dedicated to the proof of the main results. Eventually, we state %give some extensions 
in Section \ref{EXTENSIONS} %, namely we state 
weaker assumptions under which our main result holds and we also briefly discuss how to extend it to the case of a stable process perturbed by a compound Poisson process.

\mysection{Assumptions and Main results}
\label{ASS_and_RES}

\subsection{Assumptions and Notation}
 \label{ASS_1}
In the following we consider \textit{symmetric} stable processes,
that is, for all $t\ge 0 $, $u\in \R^d $,
\begin{eqnarray}
\label{TF_STABLE}
\E[\exp(i\langle u, Z_t\rangle)]
&=&\exp(it\langle \gamma,u\rangle +t\int_{S^{d-1}}^{}\bint{0}^{+\infty}\left(e^{i\rho\langle u,s\rangle}-1-i\frac{\langle u, \rho s \rangle}{1+\rho^2}\right) \frac{d\rho}{\rho^{1+\alpha}}\tilde \lambda(ds))\nonumber\\
&=&\exp(it\langle \gamma,u\rangle-t\int_{S^{d-1}}^{}|\langle s,u \rangle|^\alpha\lambda(ds)), 
\end{eqnarray}
where $\tilde \lambda $ is a symmetric measure on the unit sphere $S^{d-1} $ (i.e. for every $A$ in the Borel $\sigma$-field ${\cal
B}(S^{d-1}), \ \tilde \lambda(A)= \tilde \lambda(-A) $). The second equality in
equation \eqref{TF_STABLE} is then obtained by direct integration over $\rho$ and $\lambda=C_\alpha \tilde \lambda $ with 
$$C_\alpha:=%\left( \Gamma(1-\alpha)\alpha^{-1}\I_{\alpha\in(0,1)}+\Gamma(2-\alpha)(\alpha(\alpha-1))^{-1}\I_{\alpha\in(1,2)} \right)|\cos\left(\frac{\pi\alpha}{2}\right)|+\frac{\pi}{2}\I_{\alpha=1}.
\Gamma(1-\alpha)\alpha^{-1}\cos\left(\frac{\pi\alpha}{2}\right)\I_{\alpha \neq 1}+\frac{\pi}{2}\I_{\alpha=1}.$$
We refer to the proof of Theorem 9.32 in Breiman \cite{brei:68}  and 
Lemma 2, Chapter XVII.4 in Feller \cite{fell:66:2} for the expression of $C_\alpha$.  
 
 %From It\^o's formula for semimartingales, the canonical decomposition of $Z$ (see Theorem 2.43 Chapter II in Jacod and Shiryaev
%\cite{jaco:shir:87}) and the symmetry of $\lambda$ one derives that 
 
%%%%%% On regroupe les hypotheses comme indiquŽ par le referee #2.
We now introduce our assumptions.  Fix an integer $q\ge 2 $. We assume that 
\begin{trivlist}
\item{\A{A-1}} 
For $d\ge 2$, the spherical measure $\lambda$ has a $C^q(S^{d-1})$ surface density and for all $d\ge 1 $,
 there exist constants $0< C_1\le C_2<+\infty,\ \forall p\in \R^d,$
 $$ C_1 |p|^\alpha\leq
\int_{S^{d-1}}\left|\langle p,s\rangle \right| ^{\alpha
}\lambda \left( ds\right) \leq C_2|p|^\alpha.$$ 
%with $\lambda $ defined after \eqref{TF_STABLE}.
%Note that for $d=1$, with the convention $S^0=\{-1,1\} $, the symmetry gives $\lambda(A)=C_\lambda(\delta_{-1}+\delta_1)(A), \forall A\in {\cal B}(\R)$ with $C_\lambda   >0 $. Equation \eqref{TF_STABLE} yields
%$$\forall u\in\R, \widehat \mu(u)=\E[\exp(i uZ_1)]=\exp(i u\gamma -2C_\lambda |u|^\alpha).$$
%The  second point of the above assumption thus holds with $C_1=C_2=2C_\lambda $.

\item{\A{A-2}}  The coefficients $b$ and $f$ and their derivatives up to order $q$ are uniformly bounded in $x$. Thus, for $1<\alpha <2$, $B (
x):=b(x)+f(x)\gamma $ is uniformly bounded. We impose for $0<\alpha \leq 1$,
$B( x) =0$ for all $x\in\R^d$. 
\end{trivlist}

\begin{trivlist}
\item{\A{A-3}}
There exist constants $0<\underline{c}\le \overline{c} <+\infty$ s.t. for all $x\in\R^d,\xi\in\R^d$, $$\underline{c}|\xi|^2\le\langle f(x)\xi,\xi\rangle \le\overline{c}|\xi |^2.$$
\end{trivlist}

From now on we assume that Assumptions \A{A-1}-\A{A-3} are in force.

\begin{REM}
Note that for $d=1$, with the convention $S^0=\{-1,1\} $, we have $C_1=C_2$ in \A{A-1} even without symmetry. The symmetry is actually not needed in that case, see the beginning of Section 3 in \cite{kolo:97}.
\end{REM}
\begin{REM}
The zero drift condition in \A{A-2} comes from the fact that for $\alpha \in (0,1]$ the addition of a drift of order $t$ does not correspond to a negligible term in small time with respect to the natural scale $t^{1/\alpha}$, see Appendix B in \cite{kona:meno:10:elec} for details.
\end{REM}

In the following we denote by $C$ a positive generic constant that can depend on $\alpha,d$, the bounds appearing in the previous assumptions but neither on $N$ nor on the spatial points involved. Its value may change from line to line.
Other possible dependencies, especially w.r.t. the final time $T$ are explicitely specified. 
%%%%%%%% 04/05/09: on rajoute explicitement les definitions des espaces fonctionnels.
Concerning functional spaces, we denote %by $C_b(\R^d)$ the Banach space of continuous bounded functions on $\R^d$ endowed with the sup-norm, by $C_{b,0}(\R^d)$ its closed subspace of functions vanishing at infinity and 
by $C_b^k(\R^d),\ k\in \N^*$, the Banach space of continuous bounded functions having bounded derivatives up to and including the order $k$ with the norm $\|f\|:=\max_{0\le l\le k} \sup_{x\in\R^d}|f^{(l)}(x)| $.  Eventually $C_0^k(\R^d) $ stands for the functions in $C_b^k(\R^d)$ with compact support. 
%%%%%%%%% Il reste a checker la dependance des constantes que l'on fait apparaitre, en effet celles de l'appendice $B$ sont independantes du temps. Celles des principaux theoremes non.

\subsection{Generator  }
\label{SEC_SYMB_GEN}
From equation \eqref{TF_STABLE} and standard computations, see e.g. equation (5.11) in \cite{jaco:etal:05}, we derive that for every smooth function $g \in C_0^2(\R^d)$, the generator of \eqref{GEN_EDS} writes 
\begin{eqnarray*}
\Phi g(x)=\langle B(x),\nabla_x g(x)\rangle -\bint{\R^d}^{} g(x+f(x)y)-g(x)-\frac{\langle \nabla_x g(x), f(x) y\rangle}{1+|y|^2} \nu (dy),
\end{eqnarray*}
 where $B(x)=b(x)+f(x)\gamma $ and $\nu$ stands for the L\'evy measure of $Z$. 
 Introduce for all $A\in {\cal B}(\R^d),\ \nu_{f(x)}(A):=\nu(\{y\in\R^d: f(x)y\in A \}) $ and denote by $\tilde \lambda_{f(x)}$ its spherical part (which is still a symmetric measure). Setting $z=f(x)y$ in the above equation, using the symmetry and the polar coordinates we derive:
\begin{eqnarray}
\label{generator}
\Phi g(x)= \langle B(x),\nabla_x g(x)\rangle +\nonumber\\
\bint{S^{d-1}}^{}\int_{0}^{+\infty}\left( g(x+\rho s)-g(x) -\frac{\rho \langle \nabla_x g(x),s\rangle}{1+\rho^2}  \right)  \frac{d\rho}{\rho^{1+\alpha}} \tilde \lambda_{f(x)}(ds).%,\nonumber\\
\end{eqnarray}
%where $\tilde \lambda_{f(x)}$ stands for the spherical part of $\nu_{f(x)}(.)$. 
%In the sequel we will also need the expression of the characteristic function of 
%\begin{eqnarray*}
%\bint{\R^d}^{}  \tilde \lambda_{f(x)}\frac{d\rho}{\rho^{1+\alpha}}=
%\end{eqnarray*}
%where similarly to \eqref{TF_STABLE}, $\lambda_{f(x)}=C_\alpha \tilde \lambda_{f(x)} $.
\begin{REM}
Denote similarly to \eqref{TF_STABLE}, $\lambda_{f(x)}=C_\alpha \tilde \lambda_{f(x)} $.
The uniform ellipticity condition \A{A-3} allows to have good controls on the measure 
$\lambda_{f(x)}(\cdot)$. As a consequence of \A{A-1}, \A{A-3} one gets that  
there exist constants $0<\underline{C}_1=\underline{C}_1(\underline{c},d,\alpha)\le \overline{C}_2=\overline{C}_2(\overline{c},d,\alpha) <+\infty$ s.t. $\forall p\in\R^{d}$, $x\in \R^{d}$, 
\begin{equation}
\label{CD_ND_MES_SP}
%0<
\underline {C}_1|p|^\alpha\leq
\int_{S^{d-1}}\left| \langle p,s\rangle \right| ^{\alpha
}\lambda_{f(x)} \left(ds\right) \leq \overline{C}_2|p|^\alpha.
\end{equation}
%where analogously to \eqref{TF_STABLE} $\lambda_{f(x)}=C_\alpha \widetilde \lambda_{f(x)}  $.
\end{REM}

\subsection{Main results}
\begin{PROP}
\label{EX_DENSITY}
For every $t>0$ the solution $X_t$ (resp. $X_t^N $) of \eqref{GEN_EDS} (resp. \eqref{EUL_SCH}) has a density $p(t,x,\cdot)$(resp. $p^N(t,x,\cdot) $) w.r.t. the Lebesgue measure. Additionally, as a function of the space variables the density $p$ is in $C_b^q(\R^d\times \R^d) $ if $\alpha>1 $ and in $C_b^{q-1}(\R^d\times \R^d)  $ if $\alpha\le 1$. 
\end{PROP}
To state the theorem we first need some notation. Introduce for all $\xi\in \R^d $ and all smooth function $\varphi(t,x,y)$ the integro-differential operators:
%\begin{eqnarray}
%\Phi \varphi( t,x,y) =\langle B (x),\nabla_x \varphi(t,x,y)\rangle -\label{DEF_PHI}
%\int_{S^{d-1}}|\langle s,\nabla_x\rangle|^\alpha \varphi(t,x,y)\lambda_{f(x)}
%\left( d s\right),
%\end{eqnarray}
%and $ \forall \xi\in \R^d$,
\begin{eqnarray}
 \widetilde{\Phi }_{\xi }\varphi(t,x,y)=\langle B (\xi ),\nabla_x \varphi(t,x,y)\rangle +\label{DEF_PHI_TILDE}\nonumber\\
\bint{S^{d-1}}^{}\int_{0}^{+\infty}\left( g(x+\rho s)-g(x) -\frac{\rho \langle \nabla_x g(x),s\rangle}{1+\rho^2}  \right)  \frac{d\rho}{\rho^{1+\alpha}} \tilde \lambda_{f(\xi)}(ds) .  
\end{eqnarray}
With this definition we write for given $(x,y) \in \R^d $:
\begin{equation}
\label{DEF_OP_STAR}
\widetilde{\Phi }^{\ast }\varphi( t,x,y) =\widetilde \Phi_y \varphi(t,x,y),\ \forall m\in\N^*,\ \left( \widetilde{\Phi }_{\ast }\right) ^{m}\varphi(t,x,y)=\left( \widetilde{\Phi }%
_{\xi }\right) ^{m}\varphi(t,x,y)\mid _{\xi =x},
\end{equation}
 Note that we have $%\widetilde \Phi^* \varphi(t,x,y)=\widetilde \Phi_y \varphi(t,x,y),\ 
\widetilde{\Phi }_{\ast }\varphi(t,x,y)=\Phi \varphi(t,x,y)$ defined in \eqref{generator} but in general,  for $m\geq 2$,  $
\left( \widetilde{\Phi }_{\ast }\right) ^{m}\varphi(t,x,y)\neq \left( \Phi
\right) ^{m}\varphi(t,x,y)$.

%%%%%% Ceci est a garder mais pas ici-> En fait si car les noyaux apparaissent dans le theoreme
Define now, for $t>0$, the kernel 
\begin{equation}
\label{DEF_KERNEL}
H(t,x,y):=(\Phi-\widetilde \Phi_{y})\widetilde p^y(t,x,y)
\end{equation}
 where $\widetilde p^y(t,x,y) $ denotes the density at point $y$ of $\widetilde X_t=x+b(y)t+f(y)Z_t$. Note that the variable $y$ acts here twice: as the argument of the density and as a defining quantity of the process $\widetilde X_t (\equiv \widetilde X_{t,x,y}) $, i.e. the coefficients are frozen in $y$. 
Eventually we introduce the continuous and discrete convolution operators
\begin{eqnarray*}
\varphi\otimes \psi(t,x,y)&=&\int_{0}^{t}du\bint{}^{}dz \varphi(u,x,z)\psi(t-u,z,y), \forall t\in[0,T],\\
\varphi\otimes_N \psi(t,x,y)&=&\int_{0}^{t}du\bint{}^{}dz \varphi(\phi(u),x,z)\psi(t-\phi(u),z,y), \forall t\in \{(t_i)_{i\in\leftB 1,N\rightB} \},
\end{eqnarray*}
with $\phi(u) $ is defined just before \eqref{EUL_SCH} and denotes the largest discretization time lower or equal to $u$.
Also $\varphi\otimes H^{(0)}=\varphi $ and $\varphi\otimes H^{(r)}=\left(\varphi \otimes H^{(r-1)} \right)\otimes H$ stands for the $r$-fold convolution. 

\begin{THM}
\label{MAIN_RESULTS} 
Suppose  $q>d+4$.
Take $0<M\le q-(d+4)$. There exists a function $
R_M(T,x,y)$ with $\left| R_M(T,x,y)\right| \leq C_M(T)\left(\frac{1}{1+|y-x |^{d+\alpha}} \right):=\rho_{\alpha,M}(T,y-x)$ for some positive constant $C_M(T)$ such that
\begin{eqnarray*}
(p -p^N)( T,x,y) =
\sum_{l=1}^{M-1}\frac{h^l}{(l+1)!}\left[ p\otimes_N\left(
\Phi -\widetilde{\Phi }^{\ast }\right) ^{l+1}p^{d}\right] (T,x,y) -\\
-\sum_{k=1}^{M-1}\frac{h^k}{( k+1)! }\left[ p^{d}\otimes
_{N}\left( \widetilde{\Phi }_{\ast }-\widetilde{\Phi }^{\ast
}\right) ^{k+1}p^{N}\right] \left( T,x,y\right) +h^M R_M\left(T,
x,y\right)
\end{eqnarray*}
%where we use the convention that 
with $\sum_{l=1}^0\cdots=0$ and $\forall t\in\{ (t_i)_{i\in\leftB 1,N\rightB}\},\
p^{d}(t,x,y) :=\sum_{r=0}^{\infty }\left( \widetilde p \otimes_{N}H^{( r) }\right)(
t,x,y) $.
%\]
It holds that
\begin{eqnarray*}
\sum_{l=1}^{M-1}\left| \left( p\otimes _{N}\left( \Phi -\widetilde{%
\Phi }^{\ast }\right) ^{l+1}p^{d}\right) (
T,x,y) \right| &\leq& \rho_{\alpha,M}(T,y-x),\\
\sum_{k=1}^{M-1}\left| \left( p^{d}\otimes_{N}\left( \widetilde{%
\Phi }_{\ast }-\widetilde{\Phi }^{\ast }\right) ^{k+1}p^{N}\right)
( T,x,y) \right| &\leq &\rho_{\alpha,M} (T,y-x).
\end{eqnarray*}
\end{THM}
%%%%%%% Notons qu'a priori seul le controle en \rho_{\alpha\wedge\theta}(x-y) est utle dans la mesure ou il s'agit d'un resultat asymptotique, i.e. lorsque |x-y| est grand!!!
\begin{REM}
In the above expression, one writes for all $l\in\leftB 1,M-1\rightB $,
\begin{eqnarray*}
\left(
\Phi -\widetilde{\Phi }^{\ast }\right) ^{l+1}\varphi(t,x,y)=\bsum{k=1}^{l+1}C_{l+1}^k\Phi^k(-\widetilde \Phi^*)^{l+1-k} \varphi(t,x,y),
\end{eqnarray*}
whereas, $\forall k\in\leftB 1,M-1\rightB$,
\begin{eqnarray*}
\left( \widetilde{\Phi }_{\ast }-\widetilde{\Phi }^{\ast
}\right) ^{k+1}\varphi(t,x,y)&=& \bigl[\underbrace{(\widetilde \Phi_\xi-\widetilde \Phi_y)\cdots(\widetilde \Phi_\xi-\widetilde \Phi_y) }_{(k+1)\ {\rm times}}\bigr]\left.\varphi(t,x,y)\right|_{\xi=x}\\
&
=
&
\left. (\widetilde \Phi_\xi-\widetilde \Phi_y)^{k+1}\varphi(t,x,y)\right|_{\xi=x}.
\end{eqnarray*}
\end{REM}
\begin{REM}
The terms in the previous expansion depend on $N$. Anyhow using iteratively the Theorem and controls on $\otimes_N-\otimes $ (see also Lemma \ref{LEMME_DISCR}) it is possible to obtain an expansion with terms independent of $N$. For small $M$ explicit formulas are thus easily derived but in all generality the terms become less transparent. For $M=2$ one gets
\begin{eqnarray*}
(p-p^N)(T,x,y)=\frac h2\left(p\otimes_N(\Phi-\widetilde \Phi^*)^2p^d -p^d\otimes_N(\widetilde\Phi_*-\widetilde\Phi^*)^2p^N\right)(T,x,y)\\
 +h^2R_2(T,x,y)
\end{eqnarray*}
\vspace*{-25pt}
\begin{eqnarray*}
\phantom{(p-p^N)(T,x,y)}
=\frac h2 \left(p\otimes(\Phi-\widetilde \Phi^*)^2p -p\otimes(\widetilde\Phi_*-\widetilde\Phi^*)^2p\right)(T,x,y)\\
+h^2\widetilde R_2(T,x,y)%\\
= \frac h2 (p \otimes (\Phi^2-\widetilde \Phi_*^2) p)(T,x,y)+h^2 \widetilde R_2(T,x,y),
\end{eqnarray*}
where $\widetilde R_2(T, x,y)\le C(T) \rho_{\alpha,2}(T,y-x)$ for some positive constant $C(T)$.

From the above expansion and the controls on the density and its derivatives, see e.g. Theorems 3.1, 3.2 and Proposition 3.1 in \cite{kolo:97} or Lemma \ref{RESTES}, we can derive the error expansion for ${\cal E}_1(x,T,N) $ for measurable functions $g$ satisfying the growth condition $\exists C>0,\ |g(x)|\le C(1+|x|^\beta), \beta<\alpha$. In particular, we do not need the smoothness assumption on $g$ required in the approach of \cite{tala:tuba:90}, \cite{prott:tala:97}.  We recall that the expansion of ${\cal E}_1(x,T,N) $ allows from a practical point of view to improve the convergence rate of the discretization error using the Romberg Richardson extrapolation. This simply consists in observing that the expansion yields $\E[g(X_T)]-(\E[2g(X_T^{2N})]-\E[g(X_T^N)]) =O(h^2)$. The associated Monte Carlo estimator, involving a refined scheme, is then used for simulations see \cite{tala:tuba:90} for details. 

Also, the expansion can be used to study the sensitivity of ${\cal E}_1(T,x,N)$ w.r.t. $x$ without any additional assumption on $g$. This is crucial for financial applications (hedging), see e.g. Guyon \cite{guyo:06} for further developments in the diffusive case.
\end{REM}
%%%%%%%%%%%% Il faut etoffer cette section, i.e la remarque sur les applications

\mysection{Stable driven equations and their Euler scheme: existence of the density and associated parametrix expansion}
\label{CONN_SDE_IDO}

\subsection{Stable driven equation}
\subsubsection{Proof of Proposition \ref{EX_DENSITY}: existence of the density for the solution of \eqref{GEN_EDS}}
%In this section, we freely use  the notations introduced at the end of Section \ref{ASS_1} for functional spaces.

For $(X_t)_{t \ge 0}$, the existence of the density derives from Proposition 3.4 in \cite{kolo:97}, where some properties 
of the fundamental solution of $\partial_t  p(t,x,y)=\Phi p(t,x,y), p(0,x,y)=\delta(y-x)$ are discussed, and a standard identification argument, see e.g.  Dynkin \cite{dynk:63},
Theorem 2.3, p. 56. The stated smoothness of the density is then a consequence of point (ii) of the same Proposition.

\begin{REM}
The existence of the density is discussed in Bichteler \textit{et
al.} \cite{bich:grav:jaco:87}, where it is proved thanks to a
Bismut-Malliavin approach. This technique requires the computation
of a tangent equation associated to the gradient flow that involves
the derivatives of the coefficients of equation \eqref{GEN_EDS}.
Thus, some additional smoothness of the coefficients is needed, see e.g.
Theorem 6.48 of the above reference. We also mention the result of Picard \cite{pica:96}, Theorem 4.1, that gives existence and smoothness of the density for L\'evy driven SDEs for very singular L\'evy measures, provided there are sufficiently small jumps. For smooth coefficients $b,f$, it includes in particular the case of \eqref{GEN_EDS} where the spherical measure $\lambda $ can be atomic.  %% Rq: a priori toujours Malliavin oriented.
\end{REM}

\subsubsection{Parametrix expansion of the density}
For the sake of completeness and also because it is crucial for the discrete model
we briefly recall how to get through a "parametrix" approach a
series expansion for the density $p(t,x,y)$. 
 
Introduce, for all $x,y\in\R^d $ the following stochastic "frozen"
stable driven equation $\widetilde{X}_t \equiv \widetilde{X}_{t,x,y}$ % Notations coherentes avec \cite{kona:mamm:00}
defined for $t\ge 0$ by
\begin{equation}
\widetilde{X}_{t}=x+\int_{0}^{t}b(y) du+\int_{0}^{t}f\left( y\right)
dZ_{u} \label{five}.
\end{equation}
By computation of the Fourier transform of $Z_t$ and Fourier inversion
the transition density $%\widetilde{p}( t,x,z)\equiv
\widetilde p^y(t,x,z)$ of $\widetilde{X}_t$ at point $z\in\R^d$ 
explicitly writes 
\begin{eqnarray}
\widetilde{p}^y( t,x,z) =\frac 1{\left( 2\pi \right)
^d}\int e^{-i\langle z-x-t B ( y) ,p\rangle }%\nonumber \\
\exp
\left\{ -t
\int_{S^{d-1}}\left| \left\langle p, s\right\rangle
\right| ^\alpha \lambda_{f(y)} \left(d s\right) \right\} dp,\label{DENS_GELEE}
\end{eqnarray}
where $\lambda_{f(y)} $ has been introduced in Section \ref{SEC_SYMB_GEN}.
The densities of the solutions of \eqref{five} and \eqref{GEN_EDS} satisfy respectively
\begin{eqnarray}
\frac{\partial \widetilde p^y}{\partial t}(t,x,z)&=&\widetilde \Phi_y \widetilde p^y(t,x,z),\text{ for }t>0,\ (x,z)\in(\R^d)^2, \ \widetilde p^y(
0,x,z) =\delta ( z-x),\nonumber  \\
 \frac{\partial p}{\partial t}(t,x,z)&=&\Phi p(t,x,z),\text{ for }t> 0,\ (x,z)\in (\R^d)^2, \ p
( 0,x,z) =\delta ( z-x).  %\nonumber\\
 \label{thirteen}
\end{eqnarray}
%Note that for the first equation above the $y $ in $\widetilde \Phi_y $ is the same as the one associated to $\widetilde p(t,x,z)\equiv \widetilde p^y(t,x,z)$.
Note carefully that the derivatives in $\widetilde \Phi_y $ are taken w.r.t. the $x$ variable.

We will speak about the operators appearing in (\ref{thirteen}) as the ''frozen'' and ''unfrozen'' ones.
%From \eqref{thirteen}, we derive a parametrix
%expansion for the density $p(t,x,\cdot)$ of the solution of
%\eqref{GEN_EDS}. 
In the following $\forall (t,x,z)\in \R^{+*}\times (\R^d)^2, \widetilde p(t,x,z):=\widetilde p^z(t,x,z) $. Hence, from \eqref{DEF_KERNEL}  $\forall (t,z,y)\in\R^{+*}\times (\R^d)^2,\ H(t,z,y)=(\Phi-\widetilde \Phi_y)\widetilde p(t,z,y)=(\widetilde \Phi_z-\widetilde \Phi_y)\widetilde p(t,z,y) $. 
\begin{PROP}[Parametrix expansion of the density]
 With the notations of Section \ref{SEC_SYM_GEN}, the following representation holds
\begin{equation}
p( t,x,y) =\sum_{r=0}^{\infty }( \widetilde{p}\otimes H^{\left( r\right) }) \left( t,x,y\right).
\label{eighteen}
\end{equation}
\end{PROP}
\textit{Proof.} %Observe that $\Phi =\widetilde \Phi_x $. From the previous definition of $H$, the function $\widetilde p $, fundamental
%solution of \eqref{thirteen}, satisfies the equation
%\begin{equation}
%\frac{\partial \widetilde{p}}{\partial t}\left(
%t,x,y\right) =\widetilde \Phi_x \widetilde{p}_{
%}(t,x,y)-H( t,x,y). \label{fifteen}
%\end{equation}
Equations \eqref{thirteen} correspond to the forward Kolmogorov equations. Consider now the backward equation for $p $, namely, $\partial_s p(s,x,z)= {}^t\widetilde \Phi_z p(s,x,z)$  where $^t \widetilde \Phi_z  $ stands for the adjoint operator of  $\widetilde \Phi_z $ and the derivatives are taken w.r.t. $z$.  Differentiating under the integral we have from (\ref{thirteen})
\begin{eqnarray*}
(p-\widetilde{p})\left(
t,x,y\right) =\int_{0}^{t}ds\frac{\partial }{\partial s}\left[ \int
p\left( s,x,z\right) \widetilde{p}\left(
t-s,z,y\right) dz\right] =\\
 \int_{0}^{t}ds\int [\left( ^{t}\widetilde \Phi_z 
p\right)(s,x,z) \widetilde{p}\left( t-s,z,y\right)
-p\left( s,x,z\right) \widetilde \Phi_y \widetilde{p}\left( t-s,z,y\right) ]dz%\\
=p\otimes
H ( t,x,y).
\end{eqnarray*}
The
representation \eqref{eighteen} then follows by simple
iteration.\finpreuve
\begin{REM}
Note that the previous expansion is "formal". The convergence of the r.h.s. in \eqref{eighteen} is investigated in the proof of Theorem 3.1 in \cite{kolo:97} and can also be derived with the controls of Lemmas \ref{CTR_DER_TILDE_U_ALPHA} and \ref{LEMME_DER_ANNEXE} below. For the sake of completeness, a short proof of this convergence is also given in Appendix \ref{CTR_SERIE_PARAM}.
\end{REM}

\subsection{Euler scheme}
We consider now, for given $N\in\N^*$, the Euler scheme for equation
(\ref{GEN_EDS}) at the discretization times:
$$
X_{0}^{N}=x,X_{t_{i+1}}^{N}=X_{t_i}^{N}+b(X_{t_i}^N)h+f\left(
X_{t_{i}}^{N}\right) \left( Z_{t_{i+1}}-Z_{t_i}\right)
$$
recalling $h=T/N$. 
\subsubsection{Proof of Proposition \ref{EX_DENSITY} for the Euler scheme: existence of the density}
For each $N\in\N^*$, $(X_{t_i}^{N})_{i\in\leftB
0,N\rightB}$ is a Markov chain. Given the past $ \{X_{t_l}^{N}=x_{l},\
l\in\leftB 0,i\rightB\}$, the conditional distribution of the
innovations $b(X_{t_i}^N)h+f\left( X_{t_i}^{N}\right) \left(
Z_{t_{i+1}}-Z_{t_i}\right)$ has conditional density $\widetilde{p}^{x_i}(h,0,\cdot)$ (with the notation of \eqref{five}, \eqref{DENS_GELEE}).
This proves the existence  of the density for the discretization scheme.
\subsubsection{Parametrix expansion for the Euler scheme}
To give for the Euler
scheme an expansion similar to equation \eqref{eighteen}, that will also be the
starting point for our error expansion, we need to define, for fixed
$j,k,\ 0\leq j<k\leq N$ and $x,y\in \R^{d} $ additional ''frozen''
Markov chains $(\widetilde{X} _{t_l}^{N})_{l\in\leftB
j,k\rightB}=(\widetilde{X}_{t_l,x,y}^{N})_{l\in\leftB j,k\rightB}$.
Their dynamics is described by
\[
\widetilde{X}_{t_j}^{N}=x,\widetilde{X}_{t_{i+1}}^{N}=\widetilde{X}%
_{t_i}^{N}+b(y)h+f(y) \left( Z_{t_{i+1}}-Z_{t_i}\right) ,\
i\in\leftB j,k-1\rightB.
\]
Given the past $\{\widetilde{X}_{t_l}^{N}=x_{l},l\in\leftB
j,i\rightB\}$, the conditional distribution of the innovations
$b(y)h+f( y)( Z_{t_{i+1}}-Z_{t_i}) $ has conditional density
$\widetilde{p}^{y}( h, 0,\cdot)   $ and, hence, does not depend on the past.
Note that for the grid points $(t_i)_{i\in\leftB 0,N\rightB} $ the
transition densities of the solution $\widetilde{X}_{s,x,y}$ of \eqref{five}
coincide with the transition densities of the chain $
\widetilde{X}_{t_j,x,y}^{N}$ for $N\in\N^*,\ x,y\in \R^{d}$ and
$s=t_j$.

For all $ 0\le j<k\le N,\ (x,y)\in (\R^d)^2 $, we denote  by $p^{N}\left( t_k-t_j,x,y\right) $ and
%\begin{eqnarray}
%\label{dens_CDM_FRO}
$\widetilde{p}^{N}(t_k-t_j,x,y)$
% =\widetilde{p}^y(t_k-t_j,x,
%y)=\widetilde p(t_k-t_j,x,y)
%\end{eqnarray}
the transition probability densities between times $t_j $ and $t_k$ from point $x$ to $y$
of the chains $X^{N}$ and $ \widetilde{X}^{N}$ respectively. In particular,
\begin{eqnarray}
\label{dens_CDM_FRO}
\widetilde{p}^{N}(t_k-t_j,x,y) =\widetilde{p}^y(t_k-t_j,x,
y)=\widetilde p(t_k-t_j,x,y).
\end{eqnarray}

Before stating the parametrix expansion of $p^N $ in terms of
$\widetilde p^N$, we need to introduce a kernel $H_{N}$ that is the
"discrete" analogue of $H$ defined in \eqref{DEF_KERNEL}:
\begin{equation}
H_{N}(t_k-t_j,x,y)=\left\{ L_{N}-\widetilde{L}_{N}^y\right\} \widetilde{p}^{N}(t_k-t_j,x,y),
\label{twenty three}
\end{equation}
with
\begin{eqnarray*}
 L_{N}\varphi(t_k-t_j,x,y) =h^{-1}\{\int p^{N}(h,
x,z) \varphi( t_k-t_{j+1},z,y) dz- \varphi(t_k-t_{j+1},x,y) \},\\
\widetilde{L}_{N}^y\varphi(t_k-t_j,x,y) = h^{-1}\{\int \widetilde{p}^{y}(h, x,z)\varphi (t_k-t_{j+1},z,y) dz-\varphi(t_k-t_{j+1},x,y) \}.
\end{eqnarray*} 
Note that the previous definitions yield $p^N(h,x,z)=\widetilde p^x(h,x,z)$. We also mention that, because of the discretisation, there is a slight "shift" in time in the definition of $H_N$. Namely we have
$t_k-t_{j+1} $ instead of the somehow expected $t_k-t_j $.

\begin{LEMME}
For $0\leq j<k\leq N$ the following formula holds:
\begin{equation}
\label{para_CDM} p^{N}( t_k-t_j,x,y)=\sum_{r=0}^{k-j}(
\widetilde{p}^{N}\otimes _{N}H_{N}^{( r) }) ( t_k-t_j ,x,y )
\end{equation}
where in the calculation of $\widetilde{p}^{N}\otimes
_{N}H_{N}^{\left( r\right) }$ ($r$-fold convolution) we define
\[
p^{N}( 0,x,y) =\widetilde{p}^{N}( 0,x,y) =\delta ( x-y) .
\]
\end{LEMME}
The proof of this lemma is given in \cite{kona:mamm:00}, Lemma 3.6
and does not rely on the specific distribution of the innovations.
\begin{REM}
With the convention that $H_N^{(r)}=0 $ for $r>k-j $, equation \eqref{para_CDM} also writes 
%\begin{equation*}
$%\label{para_CDM} 
p^{N}( t_k-t_j,x,y)=\sum_{r=0}^{\infty}(
\widetilde{p}^{N}\otimes _{N}H_{N}^{( r) }) ( t_k-t_j ,x,y )$.
%\end{equation*}
This expression will often be used in the sequel. 
\end{REM}

\mysection{Proof of the main results}
\label{PROOFS_GB}
In this section, we state in Subsection \ref{STATE} the various
points needed to prove Theorem \ref{MAIN_RESULTS}. The proofs are
postponed to Subsection \ref{PROOFS}. As mentioned earlier, the key
idea consists in comparing the parametrix expansions of the
densities $p$ and $p^N$ respectively given by \eqref{eighteen} and
\eqref{para_CDM}. \textit{In the whole section we suppose that the assumptions of
Theorem \ref{MAIN_RESULTS} hold}.

\subsection{Proof of Theorem \ref{MAIN_RESULTS}}
\label{STATE}
 For the previously mentioned comparison to be possible we first need
 to estimate a difference between the transition density $
p( T,x,y) $ and $p^{d}(T,x,y):=\bsum{r\ge 0}^{} \widetilde p\otimes_N H^{(r)}(T,x,y)$ which is the analogous of \eqref{eighteen} up to the discrete time convolution (i.e. $\otimes $ replaced by
$\otimes _{N}$). We refer to \eqref{generator}, \eqref{DEF_PHI_TILDE}, \eqref{DEF_OP_STAR}, \eqref{DEF_KERNEL} for the definition of operators and kernels. 
\begin{LEMME}[Time discretization]
\label{LEMME_DISCR}
One has:
\begin{eqnarray*}
(p-p^{d})( T,x,y)
&=&\sum_{l=1}^{M-1}\frac{h^l}{(l+1)!}\left( p\otimes
_{N}\left( \Phi -\widetilde{\Phi }^{\ast }\right) ^{l+1}p^{d}\right) \left( T,x,y\right)\\
&& + h^MR_{M,1}\left( T,x,y\right)
\end{eqnarray*}
with
\begin{eqnarray*}
\sum_{l=1}^{M-1}\left| \left( p\otimes _{N}\left( \Phi -%
\widetilde{\Phi }^{\ast }\right) ^{l+1}p^{d}\right) \left(
T,x,y\right) \right| +  \label{twenty-eight}
| R_{M,1}(
T,x,y)| \leq \rho_{\alpha,M}(T,y-x). %\widetilde{p}(T,x,y)
\end{eqnarray*}
%for some constant $C:=C(T)$.
\end{LEMME}
Then the comparison between $p^d $ and $p^N$ is controlled
with the following
\begin{LEMME}[Comparison of the discrete convolutions]
\label{COMP_DISCR_DENS}
The following expansion
holds:
\begin{eqnarray*}
( p^{d}-p^{N})( T,x,y) &=&-\sum_{k=1}^{M-1}
\frac{h^k}{\left( k+1\right)!}\left[ p^{d}\otimes
_{N}\left( \widetilde{\Phi }_{\ast }-\widetilde{\Phi }^{\ast
}\right) ^{k+1}p^{N}\right] \left( T,x,y\right)\\
&& +h^MR_{M,2}\left( T,x,y\right)
\end{eqnarray*}
where
\begin{eqnarray*}
R_{M,2}(T,x,y) &=&-\frac{1}{M!}%
\int_{0}^{1}( 1-\tau ) ^{M}\left[ p^{d}\otimes
_{N}\left( \widetilde{\Phi }_{\ast }-\widetilde{\Phi }^{\ast }\right) ^{M+1}%
\widetilde{p}_{\tau }^{\Delta }\right] ( T,x,y ) d\tau, 
\end{eqnarray*}
\vspace*{-.3cm}
\begin{eqnarray*}
\forall t\in\{(t_i)_{i\in\leftB 1,N\rightB}\} ,\ \widetilde{p}_{\tau }^{\Delta }\left( t_i,x,y\right) &=&\sum_{r=0}^{\infty }%
\widetilde{p}_{\tau}\otimes _{N}H_{N}^{\left( r\right) }\left(
t_i,x,y\right) ,\widetilde{p}_{0}^{\Delta }=p^N,\\ %\widetilde{u}_{\alpha }
%,\ 
{\rm and}\ \forall \tau\in[0,1],\ %\\   
\widetilde p_{\tau}(t,x,y)&=&\bint{\R^d}^{} \widetilde p^x(\tau h,x,z)\widetilde p^y(t-\tau h, z,y)dz. 
%%%%%% Ceci est la vieille version qui introduisait une mauvaise convolution 
%\widetilde p_{\tau}(t,x,y)%&=&\widetilde p(y-x;t\lambda(y,\cdot)+\tau h\Delta \lambda^x(y,\cdot), t\gamma(y)+\tau h \Delta \gamma^x(y) )\\
%&:=&\int_{\R^d}^{}\widetilde p^y(t,x,z) \widetilde p^{\Delta \lambda^{x,y},\Delta \gamma^{x,y}}(\tau h,z,y)dz,\\
%\widetilde p^{\Delta \lambda^{x,y},\Delta \gamma^{x,y}}(\tau h,z,y)&:=&\frac 1{\left( 2\pi \right)
%^d}\int e^{-i\langle y-z-\tau h\Delta \gamma^{x,y} ,p\rangle }\\
%&&\times\exp
%\left\{ -\tau h
%\int_{S^{d-1}}\left| \left\langle p, s\right\rangle
%\right| ^\alpha \Delta \lambda^{x,y} \left(d s\right) \right\} dp,
% \\
%\Delta \lambda^{x,y}(ds)&:=&\lambda_{f(x)}(ds)-\lambda_{f(y)}(ds), \Delta \gamma^{x,y}:=\gamma(x)-\gamma(y).
\end{eqnarray*}
In particular $\widetilde p_0(t,x,y)=\widetilde p^y(t,x,y) $.
Also, %there exists $C:=C(T)>0$ s.t.
\begin{eqnarray*}
\sum_{k=1}^{M-1}\left| \left( p^d\otimes _{N}\left( \widetilde \Phi_*-
\widetilde{\Phi }^{\ast }\right) ^{k+1}p^N\right) \left(
T,x,y\right) \right| +  %\label{twenty-eight}
| R_{M,2}(
T,x,y)| \leq \rho_{\alpha,M}(T,y-x).%\widetilde{p}(T,x,y).
\end{eqnarray*}
\end{LEMME}
Theorem \ref{MAIN_RESULTS} is then a direct consequence of Lemmas \ref{LEMME_DISCR} and \ref{COMP_DISCR_DENS}.

\subsection{Proofs of the technical Lemmas}\label{PROOFS}
\hspace*{3cm}\phantom{BOUUUUH}\\
\textit{Proof of Lemma \ref{LEMME_DISCR}.} We start from the recurrence relation for $r\in\N^*$
\begin{eqnarray*}
\widetilde{p}\otimes H^{\left( r\right) }-\widetilde{p}\otimes _{N}H^{\left( r\right) }&=&\left[ \left( \widetilde{p}\otimes H^{\left( r-1\right) }\right) \otimes H-\left( \widetilde{p}\otimes H^{\left( r-1\right) }\right) \otimes _{N}H\right]\\
&& +\left[ \left( \widetilde{p}\otimes H^{\left( r-1\right) }\right)
-\left( \widetilde{p}\otimes _{N}H^{\left( r-1\right) }\right) %
\right] \otimes _{N}H.
\end{eqnarray*}
Summing up these terms over $r\in\N^* $ and using the linearity
of $\otimes $ and $\otimes _{N}$ we get
$
p-p^{d}=p\otimes H-p\otimes
_{N}H+\left( p-p^{d}\right) \otimes _{N}H
$.
An iterative application of this identity yields
\begin{equation}
(p-p^{d})(T,x,y)=\sum_{r=0}^{\infty }\left[ p\otimes
H-p\otimes _{N}H\right] \otimes _{N}H^{\left( r\right) }(T,x,y).
\label{twenty-nine}
\end{equation}
By definition, for all $k\in\leftB 1,N\rightB $,
\begin{eqnarray}
[ p\otimes H-p\otimes _{N}H](t_k,x,y) &=&\sum_{j=0}^{k-1}\int_{t_j}^{t_{j+1}}ds\int [p\left(
s,x,z\right) H\left(t_k-s,z,y\right)\nonumber\\
&& -p\left(t_j,x,z\right) H\left(t_k-t_j
,z,y\right) ]dz. \label{thirty}
\end{eqnarray}
A Taylor expansion of the function $\theta (s,z):=p(
s,x,z) H(t_k-s,z,y) $ in the interval $[t_j,s]\subseteq [t_j,t_{j+1}] $ gives
\begin{eqnarray}
\int [\theta(s,z)-\theta(t_j,z)]dz=
\sum_{l=1}^{M-1}\frac{(s-t_j)^l}{l!} \int  \left.\partial_\tau ^{l}\theta (\tau,z )\right|_{\tau =t_j}dz+\nonumber\\
\frac{(s-t_j)^M}{(M-1)!} \int_{0}^{1}\left( 1-\delta
\right) ^{M-1}\int \partial_\tau^{M}\theta (\tau
,z)|_{ \tau =\tau _{j}(s,\delta )}dzd\delta ,  \label{thirty-one}
\end{eqnarray}
where $\tau _{j}( s,\delta) =t_j+\delta (s-t_j).$ Note now that  $
- \partial_s p( t-s,x,z)=\Phi
p( t-s,x,z),\ \partial_t p(
t-s,x,z) ={ }^{t}\Phi p(
t-s,x,z)$. Here $^{t}\Phi={}^t\widetilde \Phi_z$ is the adjoint operator of $\Phi$ where the derivatives have to be taken w.r.t. $z$. Hence, $
\Phi p( t-s,x,z) ={ }^{t}\Phi p(
t-s,x,z)$. The same identity also holds for $\widetilde p $ with $\Phi,{}^t\Phi $ respectively replaced by $\widetilde \Phi^*,{}^t\widetilde \Phi^* $.
We therefore derive
\begin{eqnarray*}
\int  \partial_\tau\theta (\tau,z )|_{ \tau =t_j}dz&=&\int
\partial_\tau \left[ p( \tau ,x,z) %
\right] |_{\tau =t_j}H(t_k -t_j ,z,y ) dz\\
&&+\int p(t_j,x,z )  \partial_\tau
\left[ H(t_k-\tau ,z,y) \right]|_{\tau =t_j}dz\\
&=&
 \int  { }^{t}\Phi_z p( t_j,x,z) \left( \Phi -%
\widetilde{\Phi }^{\ast }\right) \widetilde{p}(t_k-t_j,z,y) dz\\
&&-\int p(t_j,x,z) \left( \Phi -\widetilde{\Phi }
^{\ast }\right) \widetilde{\Phi }^{\ast }\widetilde{p}(t_k-t_j,z,y) dz
\\
&=&\int p(t_j,x,z) \left( \Phi -\widetilde{\Phi }
^{\ast }\right) ^{2}\widetilde{p}\left( t_k-t_j
,z,y\right) dz.
\end{eqnarray*}
Iterating the differentiation we get
\begin{equation}
\int \partial_\tau^{l}\theta (\tau,z )|_{ \tau
=t_j}dz=
\int p(t_j,x,z) \left( \Phi -\widetilde{\Phi }%
^{\ast }\right) ^{l+1}\widetilde{p}(t_k-t_j
,z,y) dz,  \label{thirty-two}
\end{equation}
where we recall that for two operators $A$ and $B$ we denote by $(A-B)^{k}$ the following
sum $(A-B)^{k}=\sum_{j=0}^k C_k^j A^{k-j}(-B)^j$.

Plugging (\ref{thirty-one}) and (\ref{thirty-two}) into (\ref{thirty}) we
get
\begin{eqnarray}
[ p\otimes H-p\otimes _{N}H] (t_k,x,y) &=&
\sum_{l=1}^{M-1}\frac{h^{l}}{(l+1)!}p\otimes _{N}\left( \Phi -%
\widetilde{\Phi }^{\ast }\right) ^{l+1}\widetilde{p}\left( t_k,x,y\right)\nonumber \\
&& +h^{M}\widetilde R_{M,1}(t_k,x,y)
\label{thirty-three}
\end{eqnarray}
where
\begin{eqnarray}
\widetilde R_{M,1}(t_k,x,y) =\frac{1}{(M-1)!}\sum_{j=0}^{k-1}
\int_{t_j}^{t_{j+1}}\left[ h^{-1}( s-t_j) \right]
^{M}\int_{0}^{1}\left( 1-\delta \right) ^{M-1}\times\nonumber \\
\int \partial_\tau^{M}\left[ p\left(
\tau,x,z\right) H\left( t_k-\tau,z,y\right) \right] |_{\tau =\tau
_{j}(s,\delta )}dsdzd\delta  \label{thirty-four}.
\end{eqnarray}
Plugging (\ref{thirty-three}) and (\ref{thirty-four}) into (\ref
{twenty-nine}) we get
\begin{eqnarray}
(p-p^{d})( T,x,y)
&=&\sum_{l=1}^{M-1}\frac{h^l}{(l+1)!}\times
\sum_{r=0}^{\infty }p\otimes _{N}\left( \Phi -\widetilde{\Phi }%
^{\ast }\right) ^{l+1}\widetilde{p}\otimes _{N}H^{\left( r\right)
}\left(T,x,y\right)\nonumber\\
&& +h^M R_{M,1}\left( T,x,y\right)
\label{thirty-five}
\end{eqnarray}
with $R_{M,1}\left( T,x,y\right) =\sum_{r=0}^{\infty }( \widetilde R_{M,1}\otimes
_{N}H^{( r) }) ( T,x,y)$.  

Now we apply that for a linear operator $S$ and its adjoint $^{t}S$ we have $
p\otimes _{N}S\widetilde{p}={}^{t}Sp\otimes _{N}%
\widetilde{p}$. This gives
\begin{eqnarray*}
\sum_{r=0}^{\infty }p\otimes _{N}\left( \Phi -\widetilde{\Phi }%
^{\ast }\right) ^{l+1}\widetilde{p}\otimes _{N}H^{( r)
}( T,x,y) =\\
^{t}\left[ \left( \Phi -\widetilde{\Phi }^{\ast }\right) ^{l+1}\right]
p\otimes _{N}\sum_{r=0}^{\infty }\left( \widetilde{p}\otimes _{N}H^{( r) }\right) (T,x,y) =%\\
p\otimes _{N}\left( \Phi -\widetilde{\Phi }^{\ast }\right)
^{l+1}p^{d}(T,x,y),
\end{eqnarray*}
which plugged into (\ref{thirty-five}) gives the desired expansion.
The stated bound follows by application of the estimates given
in Lemma \ref{RESTES} below. We only give the proof for the first summand, the other terms of the sum over $l$ and the remainder $R_{M,1}(T,x,y)$ can be handled in a similar way. Write
\begin{eqnarray*}
p\otimes_N(\Phi-\widetilde \Phi^*)^2p^d(T,x,y)&=&\bsum{j=0}^{N-1}h\bint{}^{}p(t_j,x,z)(\Phi-\widetilde \Phi^*)^2 p^d(T-t_j,z,y)dz\\
&:=&S_1+S_2,
\end{eqnarray*}
where in $S_1$ (resp. $S_2$) the sum is taken over $I_1:=\{j\in\leftB 0,\lfloor \frac{N-1}{2} \rfloor \rightB\} $ (resp. $I_2:=\{j\in\leftB  \lfloor \frac{N-1}2 \rfloor+1 , N-1 \rightB \} $).
For $S_1 $ (resp. $S_2$), $p^d(T-t_j,z,y) $ (resp. $p(t_j,x,z) $) is non singular. 
%One could show 
From  Lemma \ref{RESTES} below equation \eqref{CT_OP},
%see also Appendix \ref{CTR_SERIE_PARAM}, using as well the proof of Proposition 2.5 and Theorem 3.2 in \cite{kolo:97} (\textbf{inclure estimees de type Aronson de Kolokoltsov dans notre THM 3.1?}) that 
there exists $ C:=C(T)$ s.t. for all $(x,y,z)\in(\R^d)^3 $, 
\begin{eqnarray}
p(s,x,z)\le C \widetilde p^y(s,x,z)%&
,
%&
\ p^d (s,z,y)\le C \widetilde p^y(s,z,y),\ \forall s\in]0,T], \nonumber \\
|(\Phi-\widetilde \Phi^*)^2 p^d(T-t_j,z,y)| %&
\le 
%& 
C \widetilde p^y(T-t_j,z,y), \forall j\in I_1,\nonumber \\
\left| {}^t\left[ \left(\Phi- \widetilde \Phi^* \right)^2\right] p(t_j,x,z)\right| %&
\le
%& 
C \widetilde p^y(t_j,x,z), \forall j\in I_2.\label{CT_IPP}
\end{eqnarray}
  % On a explicite les definitions
%Thus, recalling that for all $s\in [0,T],\ \widetilde p(T-s,z,y)=\widetilde p^y(T-s,z,y)$, t
The semigroup property for $\tilde p^y $ yields $|S_1|+|S_2|\le C\widetilde p(T,x,y) $. One eventually checks from Proposition \ref{PROP_CTR_PREAL} that $\widetilde p(T,x,y):=\widetilde p^y(T,x,y)\le \rho_{\alpha,M}(T,y-x) $.\finpreuve\\

%%%%%%%%%%%%%%%%%%%%%%%%%%% Attention, il faut bien specifier dans ce lemme que la variable de Fourier est la seconde.
\textit{Proof of Lemma \ref{COMP_DISCR_DENS}}. Let us denote by
$\F[\psi]( z) =\int \exp (i\langle
z,p\rangle )\psi(p)dp$ the Fourier transform of a
function $\psi$.  
Introduce now for all $u,t,\ u<t,$ $u,t\in \{(t_i)_{i\in\leftB 0,N\rightB}\}$, $p\in\R^d $,
\begin{eqnarray*}
\psi(p)&=&h(L_N-\widetilde L_N^y)\widetilde p^{y}(t-u,x,p)\\ 
&=&\int{}^{}p^{N}(h,x,w)\widetilde p^{y}(t-(u+h),w,p)dw -\widetilde p^{y}(t-u,x,p).
\end{eqnarray*}
Note that in particular according to \eqref{twenty three}, $\psi(y)=hH_N(t-u,x,y)$.
Taking the characteristic functions of the densities involved in the above equation, we obtain from \eqref{DENS_GELEE} and (\ref{twenty three})
that
$$ \F[\psi](z):=G_z(1)-G_z(0)$$
with
%%%%%%%%%%%% On conserve la vieille version 
%\begin{eqnarray*}
%G_z(\delta )=\exp \biggl[i\langle x,z\rangle+i(t-u)\langle \gamma( y) ,z\rangle +i
%\delta h\langle \Delta \gamma ^{x,y} 
%,z\rangle +\\
%\int_{S^{d-1}}\int_{0}^{\infty }\left( e^{i\left\langle z,\rho s\right\rangle }-1-\frac{i\left\langle z,\rho  s\right\rangle }{1+\rho ^{2}}
%\right) \frac{d\rho }{\rho ^{1+\alpha }}\left[ (t-u)\widetilde \lambda_{f(y)} (ds) +\delta h \Delta \widetilde \lambda ^{x,y}( ds) 
%\right] \biggr],
%\end{eqnarray*}
%%%%%%%%%%%%% Et l'on utilise la mesure spectrale spherique, i.e. integration de la composante radiale dans la precedente ecriture.
\begin{eqnarray*}
G_z(\tau )=\exp \biggl[i\langle x,z\rangle+i(t-u)\langle B( y) ,z\rangle +i
\tau h\langle \Delta B ^{x,y} 
,z\rangle \\
-\int_{S^{d-1}} |\langle z, s\rangle|^\alpha \left[ (t-u) \lambda_{f(y)} (ds) +\tau h \Delta \lambda ^{x,y}( ds) 
\right] \biggr],
\end{eqnarray*}
where $\Delta B^{x,y }=B
(x)-B( y)$, $\Delta  \lambda ^{x,y }(
ds) = \lambda_{f(x)}(d s) - \lambda_{f(y)}( d s)$. Note in particular that $\forall \tau\in [0,1]$,
\begin{eqnarray}
\label{EXP_GZ}
G_z(\tau)=G_z(0)\times \exp\left(\tau h\left[ i \langle  \Delta B^{x,y}, z\rangle-\bint{S^{d-1}}^{}|\langle z,s\rangle|^\alpha\Delta\lambda^{x,y}(ds) \right] \right).\nonumber\\ 
\end{eqnarray}
A Taylor expansion yields
$\F[\psi](z)=\sum_{k=1}^{M}\frac{1}{k!}
G_z^{(k)}(0)+\frac{1}{M!}\int_{0}^{1}(1-\tau )^{M}G_z^{(M+1)}(\tau )d\tau$.
From \eqref{EXP_GZ}, one derives that for $k\in\N^*$:
\begin{eqnarray*}
\frac{1}{k!}G_z^{(k)}(0)=\frac{h^k}{k!}G_z(0)%\times\\
\left[ i\left\langle \Delta B^{x,y}
,z\right\rangle -\int_{S^{d-1}} |\langle z,  s \rangle |^\alpha 
\Delta \lambda ^{x,y}(ds) \right] ^{k}.
\end{eqnarray*}
Observe now that $G_z(0)=\F[\theta](z)$, $\theta(p):=\widetilde{p}^y(t-u,x,p) $.
Using the well-known properties of the Fourier transform one gets for all $k\in\leftB 1,M\rightB $
\begin{eqnarray*}
G_z^{(k)}(0)=h^k\F\left.\left[  \left( \widetilde \Phi_\xi -\widetilde
\Phi_y\right)^k \theta \right]\right|_{\xi=x}(z),
\end{eqnarray*}
where the operators $\widetilde \Phi_. $ are applied w.r.t. the $x$ component and the Fourier transform is applied w.r.t. the $p$ component of $\widetilde p^y(t-u,x,p) $. Also, in the above writing,  we compute the Fourier transform for an arbitrary fixed $\xi\in\R^d$ and we then put $\xi=x $.

Hence, 
\begin{eqnarray*}
\F[\psi](z)=&\sum_{k=1}^{M}\frac{1}{k!}
G_z^{(k)}(0)+\frac{1}{M!}\int_{0}^{1}(1-\tau )^{M}G_z^{(M+1)}(\tau )d\tau =
\\ 
&\sum_{k=1}^{M}\frac{h^k}{k!}\F\left.\left[ \left( \widetilde{\Phi }_{\xi
}-\widetilde{\Phi }_{y}\right) ^{k}%\widetilde{u}_{\alpha
%}(y-x;(t-s)\lambda (y,\cdot),(t-s)\gamma (y))
\theta 
\right]\right|_{\xi=x}(z) +\\
 &\frac{h^{M+1}}{M! }{\cal F}\left[\int_{0}^{1}(1-\tau )^{M}\left.\left[\left( \widetilde{%
\Phi }_{\xi}-\widetilde{\Phi }_{y }\right) ^{M+1} \theta_\tau\right]\right|_{\xi=x}  d\tau \right](z),
\end{eqnarray*}
where $\forall \tau\in[0,1],\ \theta_\tau(p):=\bint{\R^d}^{}\widetilde p^x(\tau h,x,z)\widetilde{p}^y(t-u-\tau h,z,p) dz$. %  On conserve de la sorte des notations homogenes a ce qu'il se passe pour \theta(p). 
Taking the inverse Fourier transform and putting $p=y$ in the above equation, observing that $H(t-u,x,y)=(\widetilde \Phi_*-\widetilde \Phi^*)\widetilde p^y(t-u,x,y)$, 
we obtain
\begin{eqnarray}
(H_{N}-H)( t-u,x,y) =\sum_{k=1}^{M-1}\frac{h^k}{%
(k+1)!}\left( \widetilde{\Phi }_{\ast }-\widetilde{\Phi }^{\ast
}\right) ^{k+1}\widetilde{p}(t-u,x,y)+\nonumber \\
\frac{h^M}{M!}\int_{0}^{1}(1-\tau )^{M}\left( \widetilde{\Phi }_{\ast
}-\widetilde{\Phi }^{\ast }\right) ^{M+1}\widetilde{p}_{\tau
}(t-u,x,y)
d\tau  \label{forty}.
\end{eqnarray}
Recall now that %for $t\in \{ (t_i)_{i\in\leftB 1,N\rightB}\} $%
\begin{eqnarray*}
(p^{d} -p^{N})( T,x,y)
=\sum_{r=0}^{\infty }[\left( \widetilde{p}\otimes _{N}H^{\left(
r\right) }\right)   -
\left( \widetilde{p}\otimes _{N}H_{N}^{\left( r\right) }\right)
 ]( T,x,y)
\end{eqnarray*}
where we put $\left( \widetilde{p}\otimes _{N}H_{N}^{\left(
r\right) }\right) \left( T,x,y\right) =0$ for $hr>T$. Summing over  $
r\in \N$ in the identity
\begin{eqnarray*}
( \widetilde{p}\otimes _{N}H^{( r) } -\widetilde{p}\otimes _{N}H_{N}^{(
r) }) ( T,x,y) =\\
%\end{eqnarray*}
%\begin{eqnarray*}
\left( \left( \widetilde{p}\otimes _{N}H^{( r-1)
}\right) \otimes _{N}\left( H-H_{N}\right) \right) ( T,x,y) +
\\
\left( \left( \widetilde{p}\otimes _{N}H^{( r-1) }-%
\widetilde{p}\otimes _{N}H_{N}^{( r-1) }\right) \otimes
_{N}H_{N}\right) \left( T,x,y\right)
\end{eqnarray*}
one gets
\begin{eqnarray*}
(p^{d} -p^{N})( T,x,y) =\left[p
^{d}\otimes _{N}\left( H-H_{N}\right) +
 \left( p^{d}-p^{N}\right) \otimes _{N}H_{N}\right](T,x,y).
\end{eqnarray*}
By iterative application of the last identity we obtain
\begin{eqnarray*}
(p^{d}-p^{N})(T,x,y)=\sum_{r=0}^{\infty }\left[ p^{d}\otimes
_{N}\left( H-H_{N}\right) \right] \otimes _{N}H_{N}^{\left( r\right) }(T,x,y).
\end{eqnarray*}
We get from (\ref{forty}) that for all $t\in\{ (t_i)_{i\in \leftB 1,N\rightB}\} $:
\begin{eqnarray*}
\left( p^{d}\otimes _{N}\left( H-H_{N}\right) \right) (
T,x,y) =
-\sum_{k=1}^{M-1}\frac{h^k}{( k+1)!}\left[ p^{d}\otimes _{N}\left( \widetilde{\Phi }_{\ast }-\widetilde{\Phi }^{\ast
}\right) ^{k+1}\widetilde{p}\right] \left( T,x,y\right)\\
-\frac{h^M}{M!  }\int_{0}^{1}\left( 1-\tau \right) ^{M}%
\left[ p^{d}\otimes _{N}\left( \widetilde{\Phi }_{\ast }-%
\widetilde{\Phi }^{\ast }\right) ^{M+1}\widetilde{p}_{\tau}\right] (
T,x,y) d\tau.
\end{eqnarray*}
Eventually, 
\begin{eqnarray*}
( p^{d}-p^{N})( T,x,y) &=&-\sum_{k=1}^{M-1}%
\frac{h^k}{( k+1)!  }\left[ p^{d}\otimes _{N}\left(
\widetilde{\Phi }_{\ast }-\widetilde{\Phi }^{\ast }\right) ^{k+1}p^{N}\right]
\left( T,x,y\right)\\
&&+h^M R_{M,2}\left( T,x,y\right),
\end{eqnarray*}
\begin{eqnarray*}
R_{M,2}( T,x,y) =-\frac{1}{ M! }
\int_{0}^{1}\left( 1-\tau \right) ^{M}\left[ p^{d}\otimes
_{N}\left( \widetilde{\Phi }_{\ast }-\widetilde{\Phi }^{\ast }\right) ^{M+1}%
\widetilde{p}_{\tau}^{\Delta }\right] \left( T,x,y\right) d\tau\\
\forall t\in\{(t_i)_{i\in\leftB 1,N\rightB}\},\ \widetilde{p}_{\tau}^{\Delta }( t,x,y) =\sum_{r=0}^{\infty }%
\widetilde{p}_{\tau}\otimes _{N}H_{N}^{\left( r\right) }\left( t,x,y\right)
,\widetilde{p}_{0}^{\Delta }=p^N. %\widetilde{u}_{\alpha }.
\end{eqnarray*}
%\textit{I would say $ \widetilde{u}_{0}^{\Delta }=p_N$}. 
This proves the expansion part of the Lemma. The bound follows as in the previous proof from Lemma \ref{RESTES}.
\finpreuve
\bigskip

We now state Lemma \ref{RESTES} that allows to control the rests appearing in the expansions of Lemmas \ref{LEMME_DISCR} and \ref{COMP_DISCR_DENS}. Its proof is postponed to appendix \ref{LEM43}.
\begin{LEMME}
\label{RESTES}
%Let the conditions of Lemma \ref{CTR_DER_TILDE_U_ALPHA} be satisfied. Then 
Let $q>d+4 $. For all multi-indices $a,b $ s.t. $|a|+|b|<q-(d+4)$, the
following inequalities hold:
\begin{equation}
\left| D_{y}^{a}D_{x}^{b}p^{d}(t_k,x,y)\right|+ \left| D_{y}^{a}D_{x}^{b}p^N(t_k,x,y)\right|\leq
Ct_k^{-\frac{\left| a\right| +\left| b\right| }{%
\alpha }}\widetilde{p}(t_k,x,y),k \in\leftB 1,n \rightB,
\label{fifty nine}
\end{equation}
\[
\left| D_{y}^{a}D_{x}^{b}p(t,x,y)\right| \leq Ct^{-\frac{\left|
a\right| +\left| b\right| }{\alpha }}\widetilde{p}(t,x,y),0<t\leq
T.
\]
Also, $ \exists C:=C(T)$ s.t. for all $(x,y,z)\in(\R^d)^3, s\in]0,T] $, 
\begin{eqnarray}
p(s,x,z)\le C \widetilde p^y(s,x,z)
,
\ p^d (s,z,y)\le C \widetilde p^y(s,z,y),  \nonumber \\
|(\Phi-\widetilde \Phi^*)^k p^d(s,z,y)| 
\le  
C s^{-\frac{|k|}{\alpha}}\widetilde p^y(s,z,y),\nonumber\\
\left| {}^t\left[ \left(\Phi- \widetilde \Phi^* \right)^k\right] p(s,x,z)\right|\le C s^{-\frac{|k|}{\alpha}}\widetilde p^y(s,x,z).
 \label{CT_OP}
\end{eqnarray}
\end{LEMME}

\mysection{Extensions and conclusion}
\label{EXTENSIONS}
A careful examination of the proofs in the Appendices shows that the absolute continuity of $\lambda$ w.r.t. to the Lebesgue measure of $S^{d-1} $ can be removed in \A{A-1} provided the function
\begin{eqnarray*}
\zeta(t,x,y):=\frac{1}{(2\pi)^d}\bint{\R^d}^{} (\int_{S^{d-1}}^{} |\langle p,s\rangle|^\alpha\lambda_{f(x)}(ds))\exp\left(-t\bint{S^{d-1}}^{} |\langle p,s\rangle|^\alpha\lambda_{f(y)}(ds)\right)\\
\times \exp(-i\langle p,x\rangle )dp
\end{eqnarray*}
has bounded derivatives w.r.t. $x$ up to order $q$ (see Appendix \ref{CTR_SERIE_PARAM} and the statement of Theorem 3.1 in \cite{kolo:97}).
Also up to a standard perturbative argument, similar controls on the density can be obtained when we consider \eqref{GEN_EDS} driven by $(Z_s+P_s)_{s\ge 0} $ where $(P_s)_{s\ge 0}$ is a compound Poisson process with L\'evy measure $\nu_P(dz)=f(z)dz $ and $|f(z)|\le \frac{C}{1+|z|^{d+\beta}}, \beta>0 $,  see Theorem 4.1 in \cite{kolo:97}. In that case our main results remain valid up to a modification of the remainder. Indeed, it is the smallest exponent (or equivalently the largest tail) that leads the asymptotic behavior of $p(t,x,y) $ when $ |x-y|$ is large. Thus $\rho_{\alpha,M}(T,y-x) $ has to be replaced by $\rho_{\min(\alpha,\beta),M} (T,y-x)$ in Theorem \ref{MAIN_RESULTS}.  Eventually, good controls have been obtained on $p$ for stable-like processes, i.e. when the stability index in the generator $\Phi\psi(x) $ in \eqref{generator} can depend on the spatial position $x$, i.e. $\alpha $ turns to $\alpha(x)\in [\underline \alpha, \overline \alpha]$ strictly included in $(0,2] $ (see Section 5 in \cite{kolo:97}). %%%%%% rajouter l'inclusion stricte par rapport a ce qui est ecrit ici
Anyhow the processes associated to those generators cannot be approximated by a usual Euler scheme and the previous analysis breaks down. The approximation of such processes will concern further research.

\appendix

\mysection{Proof of the controls on the derivatives of the densities (Lemma \ref{RESTES})}
\label{LEM43}
To conclude the proof it remains to prove Lemma \ref{RESTES}. 
The first step is to get bounds on partial derivatives of the transition densities $
\widetilde{p}$ and $p$. The following estimates
generalize the ones obtained in \cite{kolo:97}, Propositions 2.1-2.3.

\begin{LEMME} 
\label{CTR_DER_TILDE_U_ALPHA}
Let $q>d+4$. There exists a constant $C>1$ such that the following estimates hold
uniformly for $\alpha $ in any compact subset of the interval $(0,2)$ and
for all $0<t\leq T,x,y,z\in\R^d$ and $\left| a\right| <q-\left( d+4\right)$
\begin{eqnarray}
\left| D_{z}^{a}\widetilde{p}^y(t,x,z)\right| \leq \frac{C}{%
t^{\left| a\right| /\alpha }}\widetilde{p}^y(t,x,z)
\label{forty one},\\
\left| D_{z}^{a}\widetilde{p}^y(t,x,z)\right| \leq \frac{C}{\left|
z-B ( y) t-x\right| ^{\left| a\right| }}\widetilde{p}^y(t,x,z).\label{forty two}
\end{eqnarray}
%recalling $\widetilde p(t,x,z):=\widetilde p(t,x,z;\lambda(y,\cdot),\gamma(y)) $.
\end{LEMME}
\begin{REM}
Equation \eqref{forty one} %and \eqref{forty two} 
extends to the stable case what is widely known in the Gaussian framework. Namely, each derivation of the density in space remains homogeneous to a stable density up to a multiplicative additional singularity
of order $t^{-1/\alpha} $. 
\end{REM}
\textit{Proof.} 
From now on we assume w.l.o.g. that  $d\ge 3$, the cases $d\in\{ 1,2\}$ can be addressed more directly. To proceed with the computations, we need to specify a useful change of coordinates. Namely, for a given direction $ \zeta \in\R^{d}\backslash \{ 0\}$ introduce for $p\in\R^d $ the spherical coordinates $(\rho,\vartheta,\varphi_2,\cdots,\varphi_{d-1})$, $\rho=|p | $ with first coordinate or main axis directed along $\zeta$, that is
\begin{eqnarray}
p_1&=&\rho \cos \vartheta,\ p_2=\rho \sin \vartheta \cos
\varphi_2,\ p_3= \rho \sin \vartheta \sin \varphi_2\cos \varphi
_3,...\nonumber \\
p_{d-1}&=&\rho \sin \vartheta \sin \varphi _2...\sin \varphi
_{d-2}\cos \varphi _{d-1},\nonumber \\
p_{d}&=&{\cal \rho }\sin \vartheta \sin \varphi _{2}...\sin \varphi _{d-2}\sin
\varphi _{d-1}, \label{COORD_SPHE}
\end{eqnarray}
$\vartheta \in \left[ 0,\pi \right] ,\varphi _{i}\in \left[ 0,\pi \right]
,i \in\leftB 2,d-2\rightB,\ \varphi _{d-1}\in \lbrack 0,2\pi ]$. Consider then the coordinates 
$(v,\tau,\phi)$ where $\tau =\cos \vartheta $ and $v=\rho\left|
\zeta\right|$, with $v\in\R^+,\ \tau\in[-1,1], \ \phi=(\varphi_2,\cdots,\varphi_{d-1})\in[0,\pi]^{d-3}\times [0,2\pi] $. In the following we write $p=p(v,\tau,\phi)$ for the previous r.h.s. in \eqref{COORD_SPHE} written in these new coordinates that is
\begin{eqnarray}
p_1&=&|\zeta|^{-1}v\tau,\ p_2=|\zeta|^{-1} v(1-\tau^2)^{1/2}\cos\varphi_2,\nonumber\\
 p_3&=& |\zeta|^{-1} v(1-\tau^2)^{1/2}  \sin \varphi_2\cos \varphi_3,...\nonumber\\
p_{d-1}&=&|\zeta|^{-1}v(1-\tau^2)^{1/2}   \sin \varphi _2...\sin \varphi
_{d-2}\cos \varphi _{d-1},\nonumber \\
p_{d}&=&|\zeta|^{-1}v(1-\tau^2)^{1/2}  \sin \varphi _{2}...\sin \varphi _{d-2}\sin
\varphi _{d-1}, \label{COORD_SPHE_2}
\end{eqnarray}
and $\bar p(\tau,\phi)=p(|\zeta |,\tau,\phi)$.

Without loss of generality we suppose $B(y)=0$. The first step consists in differentiating w.r.t $z$ the inverse Fourier
transform for $\widetilde{p}^y(t,x,z)$
\begin{eqnarray}
\widetilde{p}^y(t,x,z)=
\frac{1}{\left( 2\pi \right) ^{d}}\int_{\R^{d}}\exp \left\{
-t\int_{S^{d-1}}\left| \left\langle p,s\right\rangle \right| ^{\alpha
}\lambda_{f(y)}( ds ) \right\} \exp \left( -i\left\langle
p,z-x\right\rangle \right)dp.\nonumber \\
\label{forty three}
\end{eqnarray}
For $z=x$, \eqref{CD_ND_MES_SP} and standard computations directly give estimate \eqref{forty one}. Thus, in the following we also assume $z\neq x$ and use the previous spherical coordinates $(v
,\tau ,\phi)$ derived from \eqref{COORD_SPHE} setting $\zeta=z-x$ 
as the main axis. We obtain:
$$
D_z^{a}\widetilde p^y (t,x,z)=\frac{1}{( 2\pi ) ^{d}| z-x| ^{\left| a\right|
+d}}\int_{0}^{\infty }dv\ v^{| a| +d-1}\times
$$
$$
\int_{-1}^{1}d\tau \int_{[0, \pi]^{d-3}\times [0,2\pi]}d\phi \Psi( v,\tau,|a| ) \exp
\left\{ -t\frac{v^{\alpha }}{\left| z-x\right| ^{\alpha }}%
\int_{S^{d-1}}| \langle \overline{p},s\rangle|
^{\alpha }\lambda_{f(y)} ( ds) \right\} \times
$$
\begin{equation}
\tau ^{a_{1}}\left( 1-\tau ^{2}\right) ^{\frac{\left| a\right| -a_{1}+d-3}{2}%
}h(\phi ,a)\text{ }  \label{forty five},
\end{equation}
where $\overline{p}=p/\left| p\right|$, $a=(a_{1},...,a_{d})\in\N^d$ and %$ \Psi(v,\tau,|a|)=(-1)^{|
%a|/2}\cos(v\tau)\I_{|a|\ {\rm even}}+(-1)^{(|a|+1)/2}\sin(v\tau)\I_{|a|\ {\rm odd}}$,
\begin{eqnarray*}
\Psi(v,\tau,|a|)&=&(-1)^{|
a|/2}\cos(v\tau)\I_{|a|\ {\rm even}}+(-1)^{(|a|+1)/2}\sin(v\tau)\I_{|a|\ {\rm odd}},\\
h(\phi ,a)&=&\left\{\left( \cos \varphi _{2}\right) ^{a_{2}}\left( \sin \varphi
_{2}\cos \varphi _{3}\right) ^{a_{3}}\times ...\times \left( \sin \varphi
_{2}...\sin \varphi _{d-2}\cos \varphi _{d-1}\right) ^{a_{d-1}}\right. \\
&&\left.\times \left( \sin \varphi _{2}...\sin \varphi _{d-2}\sin \varphi _{d-1}\right)
^{a_{d}} \right\}\times V(\phi), \\
V(\phi) &=&\left( \sin \varphi _{2}\right) ^{d-3}\left( \sin \varphi _{3}\right)
^{d-4}\times ...\times (\sin \varphi _{d-3})^{2}\sin \varphi _{d-2}%d\varphi
%_{2}d\varphi _{3}...d\varphi _{d-1}
.
\end{eqnarray*}

We consider, first the case $\left| z-x\right| /t^{1/\alpha }\le \overline{ C}$, for a sufficiently small positive constant $\overline{C}$. 
In this case we expand the trigonometric function $\Psi(v,\tau,|a|) $ in (\ref{forty five})
in power series and
change the variable of integration $\frac{t^{1/\alpha }v}{\left|
z-x\right| }$ to $w$ in each term. This gives for all $k\in\N$,
\begin{eqnarray}
D_{z}^{a}\widetilde{p}^y(t,x,z)&=&\frac{ C_{|a|}}{ t^{\frac{|
a| +d}{\alpha }} }\biggl\{\sum_{m=0}^{k }\frac{( -1)^{m}}{(
2m+\I_{|a|\ {\rm odd}}) !}e_{m}^{|a|}\left( \frac{| z-x| }{t^{1/\alpha }}\right)
^{2m+\I_{|a|\ {\rm odd}}}\nonumber\\
&&+R_{k+1}^{|a|}\biggr\}, \ C_{|a|}=\frac{(-1) ^{
(| a|+\I_{|a|\ {\rm odd}})/2} }{( 2\pi ) ^{d}},%\nonumber\\
 \label{forty seven}
\end{eqnarray}
where $\forall m\in\leftB 1,k\rightB $,
\begin{eqnarray}
e_{m}^{|a|}&=&\int_{0}^{\infty }dw\int_{-1}^{1}d\tau \int_{[0,\pi]^{d-3}\times [0,2\pi]}d\phi \exp
\left\{ -w^{\alpha }\int_{S^{d-1}}\left| \left\langle \overline{p}%
,s\right\rangle \right| ^{\alpha }\lambda_{f(y)}(ds) \right\}\nonumber\\
&&w^{\left| a\right| +2m+d-\I_{|a|\ {\rm even}}}
\times\tau ^{a_{1}+2m+\I_{|a|\ {\rm odd}}}\left( 1-\tau ^{2}\right) ^{\frac{\left| a\right| -a_{1}+d-3%
}{2}}h(\phi ,a),\nonumber \\
|R_{k+1}^{|a|}|&\le & \frac{|e_{k+1}^{|a|}|}{(2(k+1)+\I_{|a|\ {\rm odd}})!}\left( \frac{|z-x|}{t^{1/\alpha}}\right)^{2(k+1)+\I_{|a|\ {\rm odd}}}.\nonumber\\
 \label{CTR_R}
\end{eqnarray}
To simplify the notations we omit the dependence of the
coefficients of our expansions on the direction $\zeta=z-x$.
From \A{A-1}, \A{A-2} and \eqref{CD_ND_MES_SP} one then derives the following bound:
\begin{eqnarray}
\left| e_{m}^{|a|}\right| \leq \frac{A_{d-2}}{\alpha {\underline C}_{1}^{\frac{| a|
+2m+d+\I_{|a| \ {\rm odd}}}{\alpha }}}\Gamma \left( \frac{| a| +2m+d+\I_{|a|\ {\rm odd}}}{\alpha }\right)\nonumber \\
B\left( m+\frac{a_{1}+1+\I_{|a|\ {\rm odd}}}{2},\frac{\left| a\right| -a_{1}+d-1}{2}\right).\label{CTR_E}
\end{eqnarray}
Here $A_{d-2}$ denotes the area of the unit sphere $S^{d-2}$ and $B$ is the $\beta$-function. Note that the
modulus of each term in the expansion (\ref{forty seven}) serves as an estimate of the remainder in a finite Taylor expansion%these asymptotic representations
. From \eqref{forty seven} we have
\begin{equation}
D_{z}^{a}\widetilde{p}^y(t,x,z)=\frac{C_{|a|}}{t^{\frac{\left| a\right| +d%
}{\alpha }}}\left( e_{0}^{\left| a\right| }\left( \frac{\left| z-x\right| 
}{t^{1/\alpha }}\right) ^{\mathbf{I}_{\left| a\right| odd}}+R_1^{|a|}
\right). \label{EXPR_DZA}
\end{equation}
Recall that we are considering the case $\frac{|z-x|}{t^{1/\alpha}}\le \overline{C} $. By Proposition 3.1 (i) from \cite{kolo:97} for some $\widetilde C$ depending on $\overline{C}$,
$\widetilde C^{-1}t^{-d/\alpha }\leq \widetilde{p}^y(t,x,z)\leq
\widetilde Ct^{-d/\alpha }$. Hence, equations \eqref{EXPR_DZA}, \eqref{CTR_E}, \eqref{CTR_R} yield
 \begin{equation}
\label{CONCLU_PETIT}
\left| D_{z}^{a}\widetilde{p}^y(t,x,z)\right| \leq \frac{C}{t^{\frac{%
\left| a\right| }{\alpha }}}\widetilde{p}^y(t,x,z)\leq \frac{C%
\overline{C}^{\left| a\right| }}{\left| z-x\right| ^{\left| a\right| }}%
\widetilde{p}^y(t,x,z).
\end{equation}

To estimate $D_{z}^{a}\widetilde{p}^y(t,x,z)$ for $\left|
z-x\right| /t^{1/\alpha }\ge {(\overline{C})}^{-1}$ we proceed as in Proposition 2.3 of \cite{kolo:97}. This gives the following representation
$
D_{z}^{a}\widetilde{p}^y(t,x,z)=\left[ D_{z}^{a}\widetilde{p}^{y }(t,x,z)\right] _{1}+\left[ D_{z}^{a}\widetilde{p}^y(t,x,z)
\right] _{2}
$
with
\begin{eqnarray}
[ D_{z}^{a}\widetilde{p}^y(t,x,z)] _{j}=\frac{1}{(2\pi)^{d}}\int_{0}^{\infty
}d\rho \rho ^{| a| +d-1}\int_{-1}^{1}d\tau \Psi(\rho
| z-x|, \tau,|a|) \times
\nonumber\\
f_{j}(\tau )\int_{[0,\pi]^{d-3}\times[0,2\pi]}d\phi \exp \left\{ -t\rho ^{\alpha }g_{\lambda_f
}(\tau ,\phi ,y)\right\} h(\phi ,a),\ j=1,2,\label{forty nine}\\
g_{\lambda_f}(\tau,\phi,y):=\int_{S^{d-1}}^{} |\langle \bar p(\tau,\phi) ,s\rangle |^\alpha \lambda_{f(y)}(ds),\nonumber
\end{eqnarray}
using the notations introduced after \eqref{COORD_SPHE}.
Here
$$
f_{1}(\tau )=\tau ^{a_{1}}\left( 1-\tau ^{2}\right) ^{\frac{\left| a\right|
-a_{1}+d-3}{2}}\chi (\tau ),\ f_{2}(\tau )=\tau ^{a_{1}}\left( 1-\tau
^{2}\right) ^{\frac{\left| a\right| -a_{1}+d-3}{2}}(1-\chi (\tau ))
$$
where $\chi (\tau )$ is a $C^{\infty }$ even truncation function $\R\rightarrow \left[ 0,1%
\right] $ that equals $1$ for $\left| \tau \right| \leq 1-2\varepsilon ,$
and $0$ for $\left| \tau \right| \geq 1-\varepsilon $ for some $\varepsilon
\in (0,\frac{1}{2})$. Because of the symmetry in $\tau $, it is easy to see that
the integral in (\ref{forty nine}) is non-zero
only if $a_{1}$ and $\left| a\right| $ are both even or odd. Expanding the
exponential at order 2 in (\ref{forty nine}) and making the change of variables $\rho
\left| z-x\right| =v$ we get
\begin{equation}
\left[ D_{z}^{a}\widetilde p^y(t,x,z)\right] _{1}=\frac{C_{|a|}}{\left| z-x\right|
^{\left| a\right| +d}}\sum_{m=0}^{2}\frac{1}{m!}b_{m}^{|a|}\left( \frac{t}{\left|
z-x\right| ^{\alpha }}\right) ^{m},  \label{fifty one}
\end{equation}
where $C_{|a|}$ is defined in \eqref{forty seven} and for $m\in\leftB 0,1\rightB$,
\begin{eqnarray*}
b_{m}^{|a|}&=&(-1)^{m}\int_{0}^{\infty }F_{m}^{|a|}(v)v^{\left| a\right| +m\alpha +d-1}dv,\\
F_{m}^{|a|}(v)&=&[\I_{|a|\ {\rm even}}{\func{Re}}-\I_{|a|\ {\rm odd}}{\func{Im}} ]\left[ \int_{-\infty }^{\infty }\exp (-i v\tau )\varphi
_{m}(\tau )d\tau \right],\\
\varphi _{m}(\tau )&=&f_{1}(\tau )\int_{[0,\pi]^{d-3}\times [0,2\pi]}g_{\lambda_f }^{m}(\tau ,\phi
,y)h(\phi ,a)d\phi.
\end{eqnarray*}
and
\begin{eqnarray*}
b_{2}^{|a|}&=&2\int_{0}^{1}\left( 1-\delta \right) \int_{0}^{\infty }F_{2,\delta
}^{|a|}(v)v^{| a| +2\alpha +d-1}dvd\delta ,
\\
F_{2,\delta }^{|a|}(v)&=&[\I_{|a|\ {\rm even}}{\func{Re}}-\I_{|a|\ {\rm odd}}{\func{Im}} ]\left[ \int_{-\infty }^{\infty }\exp (-iv\tau
)\varphi _{2,\delta }(\tau )d\tau \right],
\end{eqnarray*}
\begin{eqnarray*}
\varphi _{2,\delta }(\tau )&=&f_{1}(\tau )\int_{[0,\pi]^{d-3}\times[0,2\pi]}g_{\lambda_f }^{2}(\tau
,\phi ,y)\exp \left\{ -\delta t\left( \frac{v}{\left| z-x\right| }\right)
^{\alpha }g_{\lambda_f }(\tau ,\phi ,y)\right\}\\
&&\times h(\phi ,a)d\phi .
\end{eqnarray*}
To extend the integration to $\R$ in the definition of $(F_m^{|a|}(v))_{m\in\leftB 0,1\rightB},\ F_{2,\delta}^{|a|}(v) $, we simply use that the functions $(\varphi_m)_{m\in\leftB 0,1\rightB},\ \varphi_{2,\delta} $ have compact support in $\tau$.
However, to check that the coefficients $(b_m^{|a|})_{m\in\leftB 0,2 \rightB}$ are well defined, we have to equilibrate at infinity the term in $(v^{|a|+m\alpha+d-1})_{m\in\leftB 0,2\rightB} $. This can be done computing iterated integration by parts in $\tau$
in the definition of $(F_m^{|a|}(v))_{m\in\leftB 0,1\rightB},\ F_{2,\delta}^{|a|}(v) $. Namely, %recalling that 
if $\varphi _{m}(\tau ),m=0,1,$ and $
\varphi _{2,\delta }(\tau )$ are $C^{q}$ functions of $\tau $ with
compact support and $q>|a|+4+d>|a|+2\alpha+d $, performing $q$ integrations by parts w.r.t. $\tau$ one derives that the coefficients $(b_m^{|a|})_{m\in\leftB 0,2\rightB}$ are well defined. 
Let us now check that assumption \A{A-1} implies that $\varphi _{m}(\tau ),m=0,1,$ and $
\varphi _{2,\delta }(\tau )$ are $C^{q}$ functions of $\tau $ with
compact support. Indeed, for the unit
vectors $\overline{p}(\tau +\bigtriangleup \tau ,\phi )$ and $\ \overline{p}%
(\tau ,\phi ) $,  from elementary algebra there exists
an orthogonal matrix $A:=A(\Delta \tau)$ s.t. $\overline{p}(\tau +\bigtriangleup \tau ,\phi )=A 
\overline{p}(\tau ,\phi )$. Hence, if $\lambda _{f(x)}(ds)=\Theta(x,s)ds$ where $\Theta$ has the previous smoothness one can show 
\begin{eqnarray*}
&&\lim_{\Delta \tau \rightarrow 0}\frac{g_{\lambda_f }(\tau
+\Delta \tau ,\phi ,x)-g_{\lambda_f }(\tau ,\phi ,x)}{\Delta
\tau }\\
&=&\lim_{\Delta \tau \rightarrow 0}\frac{\int_{S^{d-1}}\{\left|
\left\langle \overline{p}(\tau ,\phi ),A^{\ast }s\right\rangle \right|
^{\alpha }-\left| \left\langle \overline{p}(\tau
,\phi ),s\right\rangle \right| ^{\alpha }\}\lambda_{f(x)} (ds)  }{\Delta \tau 
}\\
&=&\int_{S^{d-1}}\left| \left\langle \overline{p}(\tau ,\phi ),s\right\rangle
\right| ^{\alpha }\lim_{\Delta \tau \rightarrow 0}\frac{%
[\Theta(x,As)-\Theta(x,s)]}{\Delta \tau } ds\\
&=&\int_{S^{d-1}}\left| \left\langle \overline{p}(\tau ,\phi
),s\right\rangle \right| ^{\alpha }\Theta_{s}^{\prime }(x,s)\beta(\tau ,\phi ,s)ds,
\end{eqnarray*}
where $\beta(\tau ,\phi ,s)$ is $C^{\infty }$ function in $\tau $ uniformly
bounded in $(\tau ,\phi ,s)$ in our region. The process can then be iterated other $q-1$ times.

Thus all coefficients $(b_{m}^{|a|})_{m\in\leftB 0,2\rightB}$ are well defined.
%%%%%%%% Je ne vois plusbien le sens de cette assertion
%and (\ref{fifty one}) actually gives an asymptotic expansion.

Next, analogously to Proposition 2.3 in \cite{kolo:97} (where the
case $\left| a\right| =0$ was considered) and with the same rotations of the integration contours for $\alpha\in (0,1],\ \alpha\in(1,2) $, we obtain for all $k\in\N $
\begin{equation}
[ D_{z}^{a}\widetilde{p}^y(t,x,z)]_{2}=\frac{C_{|a|} }{\left| z-x\right|
^{\left| a\right| +d}}\biggl\{\sum_{m=0}^{k}\frac{1}{m!}c_{m}^{|a|}\left( \frac{t}{%
\left| z-x\right| ^{\alpha }}\right) ^{m}+R_{2,k+1}^{|a|}\biggr\}, \label{fifty two}
\end{equation}
\vspace*{-.5cm}
\begin{eqnarray*}
c_{m}^{|a|}=2[\I_{|a|\ {\rm even}}{\func{Re}}-\I_{|a|\ {\rm odd}}{\func{Im}} ][\int_{1-2\varepsilon }^{1}d\tau \int_{[0,\pi]^{d-3}\times [0,2\pi]}d\phi
h(\phi ,a)(-g_{\lambda_f }(\tau ,\phi ))^{m}
\\
\times\exp (-\frac{i\pi \alpha m}{2})(-i)^{\left| a\right| +d}\tau ^{-(\alpha
m+d+\left| a\right| )}\Gamma \left( \alpha m+d+\left| a\right| \right)
f_{2}(\tau )],
\end{eqnarray*}
%\begin{eqnarray*}
and $|R_{2,k+1}^{|a|}|\le \frac{|c_{k+1}^{|a|}|}{(k+1)!}\left(\frac{t}{|x-z |^\alpha}\right)^{k+1}
$. %\end{eqnarray*}
 Note that the
coefficients $c_{m}^{|a|}$ are also well defined because $\tau $ does not approach
zero (recall that $1-\chi(\tau)\neq 0 \Leftrightarrow |\tau|>1-2\varepsilon $). Precisely $|c_m^{|a|}|\le 2A_{d-2}C_2^m(1-2\varepsilon)^{-\alpha m+d+|a|}\Gamma(\alpha m+d+|a|) $.

Now the sum of expansions (\ref{fifty one}) and (\ref{fifty two})  gives the expansion for $D_{z}^{a}
\widetilde{p}^y(t,x,z)$. Note that by construction, the first coefficient $%
b_{0}^{|a|}+c_{0}^{|a|} $ does not depend on the spectral measure $\lambda_{f(y)}(\cdot)$
and it vanishes when the spectral measure is uniform (that is $C_{1}=C_{2}=1$
in \eqref{CD_ND_MES_SP}). This can be shown by means of representations
involving Bessel and Whittaker functions and the same rotations of the
integration contours as in Proposition 2.2 of \cite{kolo:97}, see Appendix \ref{CALC_COEFF} for details. Thus, for all $k\in\N^* $, we
get a representation
\begin{equation}
D_{z}^{a}\widetilde{p}^y(t,x,z)=\frac{C_{|a|}}{| z-x|
^{| a| +d}}\biggl\{\sum_{m=1}^{k}\frac{1}{m!}d_{m}^{|a|}\left( \frac{t}{%
\left| z-x\right| ^{\alpha }}\right) ^{m}+R_{k+1}^{|a|} \biggr\},  \label{fifty five}
\end{equation}
where $d_{m}^{|a|}=b_{m}^{|a|}+c_{m}^{|a|}$ with $b_{m}=0$ for $m\geq 3$ and $|R_{k+1}^{|a|}|\le \frac{|d_{k+1}^{|a|}|}{(k+1)!}\left(\frac{t}{|x-z |^\alpha}\right)^{k+1}$. 
Now, from Proposition 3.1 (ii) in \cite{kolo:97} $d_1^0>0$. Equation \eqref{fifty five} yields
\begin{equation*}
D_{z}^{a}\widetilde{p}^y(t,x,z)=\frac{C_{|a|}d_1^0 t}{\left| z-x\right|
^{\left| a\right| +d+\alpha}}\left( \frac{d_{1}^{\left| a\right| }}{%
d_{1}^{0}}+\widetilde R_2^{|a|}\right) ,\ |\widetilde R_2^{|a|}|\le \frac{|d_{2}^{|a|}|}{2d_1^0}\frac{t}{|x-z |^\alpha}, 
\end{equation*}
\begin{equation*}
\widetilde{p}^y(t,x,z)=\frac{C_{0}}{\left| z-x\right| ^{d}}\left( 
\frac{d_{1}^{0}t}{\left| z-x\right| ^{\alpha }}+R_2^0\right)\ge\frac{C_{0}d_1^0 t}{2\left| z-x\right| ^{d+\alpha}}
\end{equation*}
for sufficiently small $\overline{C} $.
Hence, we have
\begin{eqnarray}
|D_{z}^{a}\widetilde{p}^y(t,x,z)|&\le & \frac{CC_{|a|}}{\left| z-x\right|
^{\left| a\right| }}\frac{d_{1}^{0}t}{|z-x |^{d+\alpha}}
\leq \frac{C}{\left| z-x\right| ^{\left| a\right| }}\widetilde{p}^y(t,x,z)\nonumber \\
&\leq& \frac{C\overline{C}^{\left| a\right| }}{t^{\left| a\right|
/\alpha }}\widetilde{p}^y(t,x,z), \label{MAJ_GRAND}
\end{eqnarray}
recalling that $\frac{t^{1/\alpha}}{|z-x|}\le \overline{C} $ for the last inequality.
W.l.o.g. we can assume $\overline{C}<1 $. It remains to consider the case $|x-z |/t^{1/\alpha} \in] \overline{C}, \overline{C}^{-1} [:=I(\overline{C}) $. It follows from \eqref{forty five} that $|z-x |^{d}\widetilde p^y(t,x,z)$ and $|z-x |^{d+|a|}D_z^a \widetilde p^y(t,x,z)$ are continuous functions of $|x-z |/t^{1/\alpha}$. Since the stable density is also strictly positive, we deduce that there exists $\widetilde C $ s.t. on $I(\overline{C}) $, $|D_z^a \widetilde p^y(t,x,z)|\le  \frac{\widetilde C}{|z-x |^{|a|+d}}\le \frac{C}{|z-x|^{|a|}} \widetilde p^y(t,x,z)\le \frac{C\overline{C}^{|a|}}{t^{|a|/\alpha}}\widetilde p^y(t,x,z)  $
which concludes the proof.
\finpreuve

\begin{LEMME}
\label{LEMME_DER_ANNEXE}
%Let the conditions of Lemma \ref{CTR_DER_TILDE_U_ALPHA} be satisfied. 
Let $q>d+4$. 
There exists a 
constant $C>1$ s.t. the following estimates hold uniformly for $\alpha$ in any compact subset of the interval $(0,2) $ and for all $0<t\le T,\ x,y,v\in\R^d $ and $\left| a\right| +\left| b\right| < q-(d+4)%\wedge
%q
$:
\begin{eqnarray}
\left| D_{y}^{a}D_{x}^{b}H(t,x,y)\right| &\leq & \frac{C}{t^{\frac{\left|
a\right| +\left| b\right| }{\alpha }}}\widetilde{p}(t,x,y)\left( 1+%
\frac{\min (1,| y-x| }{t}\right), \label{fifty six}\\
\left| D_{x}^{b}H(t,x,x+v)\right| &\leq& C\widetilde{p}(t,x,x+v)\left( 1+\frac{\min (1,| v| }{t}\right) 
\label{fifty seven},\\
%\end{equation}
%\begin{equation}
\left| D_{y}^{a}D_{x}^{b}\widetilde{p}(t,x,y)\right| &\leq& \frac{C}{%
| y-B(y)t-x| ^{| a| +| b| }}%
\widetilde{p}(t,x,y)\label{fifty eight}.
\end{eqnarray}
\end{LEMME}
\textit{Proof.} Inequalities (\ref{fifty six}) and (\ref{fifty seven}) follow from the
representation
\[
H(t,x,y)=\left\langle B (x)-B (y),\nabla_x \widetilde{p}(t,x,y)%
\right\rangle +\frac{1}{(2\pi )^{d}}\int_{\R^{d}}\left| p\right|
^{\alpha }\int_{S^{d-1}}\left| \langle\overline{p},s\rangle \right| ^{\alpha
}\times
\]
\[
( \lambda_{f(y)}(ds) -\lambda_{f(x)} (ds))\exp \left\{ -t\left| p\right|
^{\alpha }\int_{S^{d-1}}\left| \langle \overline{p},s\rangle \right| ^{\alpha
}\lambda_{f(y)}(ds)\right\} \times
\]
\begin{equation}
\label{REPR_H}
\exp \left\{ -i\langle p,y-B (y)t-x\rangle \right\} dp
\end{equation}
analogously to the proof of Proposition 2.3 in \cite{kolo:97}, see also Appendix \ref{CTR_SERIE_PARAM} where \eqref{fifty six} is proved for $|a|=|b|=0 $.
Inequality (\ref{fifty seven}) contains in (3.23') p.748 of that reference.
Inequality \eqref{fifty eight} can be derived following the proof of Lemma \ref{CTR_DER_TILDE_U_ALPHA}.
\finpreuve

The proof of Lemma \ref{RESTES} can then be achieved from Lemmas \ref{CTR_DER_TILDE_U_ALPHA} and \ref{LEMME_DER_ANNEXE} adapting the arguments in Appendix \ref{CTR_SERIE_PARAM} concerning the control in terms of the frozen density for the "formal" series appearing in \eqref{eighteen}. See also the proof of Theorem 2.3 in \cite{kona:mamm:02} or Theorem 3.1 in \cite{kolo:97}.

%\appendix
\mysection{Control of the parametrix series of the density}
\label{CTR_SERIE_PARAM}
For the sake of completeness we provide in this section a complete proof of the control for the r.h.s of \eqref{eighteen} under our standing Assumptions \A{A-1}-\A{A-3}.

We first sum up in Proposition \ref{PROP_CTR_PREAL} the various estimates needed to control the convergence of \eqref{eighteen} following the proof of Theorem 3.1 in \cite{kolo:97}, namely Proposition 3.1 and its corollary, Lemma 3.1 and Propositions 3.2-3.3 of that reference. These estimates can also be directly derived from the computations of Appendix \ref{CALC_COEFF}.

\begin{PROP}
\label{PROP_CTR_PREAL}
For  all $K$ sufficiently large, there exists $C>0$ s.t. the following estimates hold uniformly for $\alpha $ in any compact subset of $(0,2)$, for all $x,y,z\in \R^d $ and for all $t\in(0,T]$.
\begin{eqnarray*}
C^{-1}t^{-d/\alpha}&\le &\widetilde p^y(t,x,z) \le Ct^{-d/\alpha},\ |x-z|\le Kt^{1/\alpha},\\
\frac{C^{-1}t}{|x-z|^{d+\alpha}}&\le &\widetilde p^y(t,x,z)\le \frac{Ct}{|x-z|^{d+\alpha}},\ |x-z|\ge Kt^{1/\alpha},\\
\widetilde p^z(t,x,z)&\le &C\widetilde p^y(t,x,z).
\end{eqnarray*}
Also, there exists $C>0$ s.t. $\forall (t,x,y)\in[0,T]\times(\R^d)^2,$
\begin{eqnarray*}
\int{}^{}dz \min(1,|z|)\widetilde p^y(t,0,z )\le Ct^\omega,\ \omega:=\min(1,1/\alpha).
\end{eqnarray*}
For all $s\in (0,t)$
\begin{eqnarray*}
\int{}^{}dz\widetilde p^z(t-s,x,z)\min(1,|y-z|)s^{-1}\widetilde p^y(s,z,y )\\
\le C\bigl(t^{-1}\min(1,|y-x |) +s^{\omega-1} \bigr)\widetilde p^y(t,x,y ),\\
%\end{eqnarray*}
%\begin{eqnarray*}
\int{}^{}dz \min(1,|z-x |)\widetilde p^z(t-s,x,z )\\
\times \min(1,|y-z|)s^{-1}\widetilde p^y(s,z,y)
\le Cs^{\omega-1}\widetilde p^y(t,x,y ),\\
\int{}^{}dz \widetilde p^z(t-s,x,z )\widetilde p^y(s,z,y )\\
\le C \widetilde p^y(t,x,y ).
\end{eqnarray*}
\end{PROP}
Introduce now for a given bounded measure $\eta$ on $S^{d-1}$ the function
%%%%%% Verifier que la notation \bar p a bien ete introduite
\begin{eqnarray*}
\varphi_\eta(t,z,\lambda_{f(y)}):=\frac{1}{(2\pi)^d}\bint{\R^d}^{}dp |p|^\alpha\bint{S^{d-1}}^{} |\langle \bar p,s\rangle|^\alpha \eta(ds)\\
\times  \exp\left(-t |p|^\alpha\bint{S^{d-1}}^{} |\langle \bar p,s\rangle|^\alpha \lambda_{f(y)}(ds)\right)\exp(-i\langle p,z\rangle).
\end{eqnarray*}
With this notation and \eqref{REPR_H} we get
\begin{eqnarray*}
H(t,x,y)=\langle B(x)-B(y),\nabla_x\widetilde p(t,x,y)\rangle +(\varphi_{\lambda_{f(y)}}-\varphi_{\lambda_{f(x)}})(t,y-x-B(y)t ,\lambda_{f(y)}).
\end{eqnarray*}
Under our standing assumptions, the mean value theorem yields $|\varphi_{\lambda_{f(y)}}-\varphi_{\lambda_{f(x)}})(t,y-x-B(y)t ,\lambda_{f(y)})|\le\min(1,|y-x|)\varphi_{\eta_{x,y}}(t,y-x-B(y)t ,\lambda_{f(y)}) $, where $\eta_{x,y} $ is a bounded measure. 
Now Proposition 2.5 in \cite{kolo:97} states that for a bounded measure $\eta$,
$$\varphi_\eta(t,z,\lambda_{f(y)})\le Ct^{-1}\widetilde p^y(t,0,z).$$
From Lemmas \ref{CTR_DER_TILDE_U_ALPHA}, \ref{LEMME_DER_ANNEXE} and the above controls one deduces $|H(t,x,y)|\le C\widetilde p^y(t,x,y)\bigl(1+t^{-1}\min(1,|x-y|)\bigr):=Cv(t,x,y)$ (which actually gives \eqref{fifty six} for $a=b=0$).

Introduce now $\beta\circ \psi (t,s,x,y)=\bint{\R^d}^{} \beta(t-s,x,z)\psi(s,z,y)dz $, i.e.
$\circ$ is the spatial part of the convolution operator $\otimes$, and set $%\hat p(t,x,y):=t^{-1}\min(1,|x-y |)\widetilde p(t,y-x;\lambda(y),\gamma(y)),%\ %\check p(t,x,y):=t\hat p(t,x,y)=\min(1,|x-y |)\widetilde p(t,y-x;\lambda(y),\gamma(y)),\ v(t,x,y):=(\widetilde p+\hat p)(t,x,y),
\widetilde v(t,x,y):=tv(t,x,y)$. From Propositon \ref{PROP_CTR_PREAL} one derives
\begin{eqnarray*}
\widetilde p\circ v(t,s,x,y)\le C(v(t,x,y)+s^{\omega-1}\widetilde p(t,x,y))\\
\widetilde v\circ v(t,s,x,y)\le C(\widetilde v(t,x,y)+(s^{\omega-1}+(t-s)^\omega )\widetilde p(t,x,y)).
\end{eqnarray*}
%%%%%%%%%%%%%%%%%%%%%%%%%%% On zappe cela pour donner directement les controles globaux qui peuvent etre factorses.
%For $\alpha\le 1 $ the previous expressions simplify to 
%\begin{eqnarray*}
%\widetilde p\circ v(t,s,x,y)\le C v(t,x,y), \widetilde v\circ v(t,s,x,y)\le C(\widetilde v+\widetilde p)(t,x,y).
%\end{eqnarray*}
Recalling $|H(t,x,y)|\le Cv(t,x,y) $, integrating the above inequalities one gets:
$$|\widetilde p\otimes H|(t,x,y)\le C(\widetilde v(t,x,y)+t^\omega \widetilde p(t,x,y)), \ |\widetilde p\otimes H\otimes H|(t,x,y)\le C^2t^\omega(\widetilde p+\widetilde v)(t,x,y) .$$
An induction yields, for all $k\in\N^*$:
\begin{eqnarray*}
|\widetilde p\otimes H^{(2k)}|(t,x,y)&\le& \frac{C^{2k}t^{k\omega}}{(k!)^2}(\widetilde p+\widetilde v)(t,x,y), \\
|\widetilde p\otimes H^{(2k+1)}|(t,x,y)&\le& \frac{C^{2k+1}t^{k\omega}}{k!(k+\I_{\alpha\in(0,1]})!}(tv(t,x,y)+t^\omega\widetilde p(t,x,y)\I_{\alpha>1}),
\end{eqnarray*}
and the the required control, i.e. $p(t,x,y)\le C \widetilde p(t,x,y)$.
The controls on the derivatives can be proved in a similar way, up to suitable rearrangements of the variable of integration, see p.747 and 748 in \cite{kolo:97}.
Also the whole proof can be carried out for $p^d,p^N $.
\finpreuve
\begin{REM}
\label{REM_FIN}
To conclude, note that by arguments similar to those used to prove Proposition \ref{PROP_CTR_PREAL}, one gets
$$|H\otimes H(t,x,y)|\le Ct^{\omega-1}\widetilde p(t,x,y), $$
which turns to be a useful estimate to derive \eqref{CT_IPP} following the above proof. 
\end{REM}

\mysection{
Additional computations concerning the derivatives of the density
}
\label{CALC_COEFF}
In this section we prove that $b_{0}^{|
a| }+c_{0}^{| a| }=0$, justifying that the first index in \eqref{fifty five} is one. 
\subsection*{Odd dimensions $d$}
%We consider first the case of odd dimensions $d$. 
From the definitions in the proof of Lemma \ref{CTR_DER_TILDE_U_ALPHA} , it is enough to show
\begin{eqnarray}
\left[ \I_{|a|\ {\rm even}}\text{Re}-\I_{|
a|\  {\rm odd}}\text{Im}\right]\left\{ \int_{0}^{\infty }\left[ \int_{\R
}^{}\exp \left( -iv\tau \right) f_{1}(\tau )d\tau \right] v^{|
a| +d-1}dv\right.\nonumber \\
\left.+2\int_{1-2\varepsilon }^{1}\left( i\tau \right)
^{-(d+\left| a\right| )}\Gamma \left( d+\left| a\right| \right) f_{2}(\tau
)d\tau\right\} =0. \label{C_TO_SHOW_1}
\end{eqnarray}
Note that since $d$ is odd, if $|a|$ is odd $i^{-(d+|a|)}=(-1)^{\frac{d+|a|}2} $ and if $|a| $ is even $i^{-(d+|a|)}=i^{-1}(-1)^{\frac{d-1+|a|}2} $. 
Hence, the contribution of the second term in \eqref{C_TO_SHOW_1} vanishes and the condition %\eqref{C_TO_SHOW} thus 
writes
\begin{equation}
\left[ \I_{| a|\  {\rm even}}\text{Re}-\I_{|
a|\ {\rm odd}}\text{Im}\right]
\left\{ 
\int_{0}^{\infty }\left[ \int_{\R}^{}\exp \left( -iv\tau \right) f_{1}(\tau )d\tau%+2\int_{1-2\varepsilon }^{1}\exp (-i\tau v)f_{2}(\tau )d\tau 
\right]
v^{| a| +d-1}dv\right\}=0.\label{C_TO_SHOW}
\end{equation}
 Denote 
$G_{1}(v)=\int_{\R }^{}\exp \left( -iv\tau \right) f_{1}(\tau
)d\tau $. 
Remind that $| a| $ and $a_{1}$ have the same parity, see proof of Lemma \ref{CTR_DER_TILDE_U_ALPHA}. % and that $d$ is odd.
\begin{trivlist}
\item[\textbf{a)}] For even $|a| ,a_{1}$, $G_{1}$ is even and belongs to a Schwartz space of
functions. Since $d$ is odd, by inverse Fourier transform, Equation \eqref{C_TO_SHOW} reduces to 
\begin{eqnarray*}
%\bsum{j\in \{1,2\}}^{}
\int_{0}^{\infty }G_{1}(v)v^{| a| +d-1}dv=\frac 12 %\bsum{j\in\{1,2 \}}^{}
\int_{\R}^{}G_{1}(v)v^{| a| +d-1}dv
\end{eqnarray*}
\begin{equation*}
=\frac{(-i)^{\left| a\right| +d-1}(2\pi)^d}{2}%\bsum{j\in\{1,2 \}}^{}
f_{1}^{( | a|
+d-1) }(0)=0.
\end{equation*}
%Now $f_{2}^{(|a|+d-1)}(0)=0$ by definition of $f_{2}$ and 
The equality
$f_{1}^{\left( \left| a\right| +d-1\right) }(0)=0$ follows from
the Leibniz differentiation rule for the product $\tau ^{a_{1}}\times
(1-\tau ^{2})^{\frac{\left| a\right| -a_{1}+d-3}{2}\text{ }}$ 
%(for this case $(1-\tau ^{2})^{\frac{\left| a\right| -a_{1}+d-3}{2}\text{ }}$ is a
%polynomial of the order $\left| a\right| -a_{1}+d-3$) 
and the identity $| a| +d-1=a_{1}+(|a| -a_{1}+d-3)+2]$. Thus \eqref{C_TO_SHOW} holds in this case.
\item[\textbf{b)}] Analogously, for odd $|a|,a_{1}$, -${\rm Im} (G_{1}(v))$ is odd and belongs to a
Schwartz space of functions. The function $( -\text{Im}G_{1}(v) v^{| a|
+d-1})$ is even. Thus 
$$%\sum_{j\in\{ 1,2\}}^{}
\int_{0}^{\infty }( -\text{Im}G_{1}(v)) v^{|a|
+d-1}dv=\frac{(-i)^{|a| +d}}{2}%\sum_{j\in\{1,2 \}}^{}
f_{1}^{( | a|
+d-1) }(0)=0$$ 
for the same previous reasons and
equation \eqref{C_TO_SHOW} holds in this case as well.
\end{trivlist}

\subsection*{Even dimensions $d$} We assume in this section that the spectral
measure is uniform with $C_1=C_2=1 $ in \A{A-1}. For $\left| a\right| $ and $a_{1}=2m$ even, equation \eqref{forty five} can be rewritten as 
\begin{eqnarray}
D_{z}^{a}\widetilde{p}^y(t,x,z)=\frac{(-1)^{| a|
/2}A_{d-2}^{a}}{( 2\pi ) ^{d}| z-x| ^{|
a| +d}}
 \sum_{j=0}^{m}(-1)^{2m-j}C_{m}^{j}\label{ANC_1A}\\
\times \int_{0}^{\infty }dvv^{|
a| +d-1}\exp \left( -t\frac{v^{\alpha }}{| z-x|^\alpha}\right) \int_{-1}^{1}(1-\tau ^{2})^{N_{j}-1/2}\cos (v\tau )d\tau
\nonumber
\end{eqnarray}
where $A_{d-2}^{a}=\int_{[0,\pi]^{d-3}\times[0,2\pi] }h(\phi ,a)d\phi $ and $N_{j}=
\frac{| a| -a_{1}+d-2+2j}{2},\ j\in\leftB 0,m\rightB$. Now recalling the definition of the Bessel function
$J_n(z):=\frac{(z/2)^n}{\Gamma(n+\frac 12)\sqrt \pi}\int_{-1}^{1}(1-t^2)^{n-1/2}\cos(zt)dt $ which is well defined for $n>1/2$ on $\C\backslash (-\infty,0)$, we get
\begin{eqnarray}
D_z^a\widetilde p^y(t,x,z)&=&\frac{(-1)^{|a|/2}A_{d-2}^{a}}{(2\pi)^{d}|
z-x| ^{|a| +d}}%
\sum_{j=0}^{m}(-1)^{2m-j}C_{m}^{j}2^{N_{j}}\Gamma (N_{j}+\frac{1}{2})\sqrt{%
\pi }\nonumber \\
&&\times \int_{0}^{\infty }v^{\frac{| a| +a_{1}+d-2j}{2}}\exp \left(
-t\frac{v^{\alpha }}{\left| z-x\right|^\alpha }\right)
J_{N_{j}}(v)dv\nonumber\\
%\end{eqnarray*}
%\begin{eqnarray}
&=&\frac{(-1)^{\left| a\right| /2}A_{d-2}^{a}}{\left( 2\pi \right) ^{d}\left|
z-x\right| ^{\left| a\right| +d}}\sum_{j=0}^{m}(-1)^{2m-j}C_{m}^{j}\Gamma
(N_{j}+\frac{1}{2})2^{N_j+1/2}\label{OLD_2}
\end{eqnarray}
\begin{eqnarray*}
\times \text{Re}\int_{0}^{\infty }\exp \left( -t\frac{v^{\alpha }}{\left|
z-x\right| ^\alpha}\right) \exp \left[ \left( \frac{1}{2}%
N_{j}+\frac{1}{4}\right) \pi i\right] W_{0,N_{j}}(2iv)v^{N_j'}dv, 
\nonumber
\end{eqnarray*}
where $W_{0,n}(z)=\frac{\exp(-z/2)}{\Gamma(n+\frac 12)}\int_{0}^{\infty}[t(1+t/z)]^{n-1/2}e^{-t}dt$, $n>1/2$, $z\in \C\backslash (-\infty,0) $, is the Whittaker function and for 
$z>0 $, $$J_n(z)=2{\rm Re}\left[\frac{1}{\sqrt{2\pi z}}\exp(\frac 12(n+\frac 12)\pi i)W_{0,n}(2iz) \right] $$ (relation (2.10) from \cite{kolo:97}) and $N_j'=\frac{|a|+a_1+d-2j-1}{2},\ j\in\leftB 0,m\rightB $. For $\alpha \in (0,1]$, from Cauchy's theorem we can change the integration path in \eqref{OLD_2} to the negative imaginary
half line. Setting then $v=-i\xi $ we obtain
\begin{equation*}
D_{z}^{a }\widetilde{p}^y(t,x,z)=\frac{(-1)^{\left| a\right|
/2}A_{d-2}^{a}}{\left( 2\pi \right) ^{d}\left| z-x\right| ^{\left|
a\right| +d}}\sum_{j=0}^{m}(-1)^{2m-j}C_{m}^{j}\Gamma (N_{j}+\frac{1}{2})2^{N_j+1/2}
\end{equation*}
\begin{equation*}
\times (-1)^{j-m}\text{Re}\left[ -i\int_{0}^{\infty }\exp \left( -t\frac{\xi ^{\alpha }%
}{\left| z-x\right|^\alpha}\exp \left( -\frac{i\pi \alpha }{%
2}\right) \right) W_{0,N_{j}}(2\xi )\xi^{N_{j}'}d\xi \right] .
\end{equation*}
Recalling the definition of $W_{0,N_j}$, we conclude expanding the exponential in power series that the first term is 0.

For $\alpha \in (1,2)$, using the same arguments we can rotate the initial contour through the angle $-\pi/(2\alpha)$. Setting then $\eta=e^{i\frac{\pi}{2\alpha}}v $ we get
\begin{equation*}
D_{z}^{a }\widetilde{p}^y(t,x,z)=\frac{(-1)^{\left| a\right| /2}A_{d-2}^{a}}{\left( 2\pi \right) ^{d}\left|
z-x\right| ^{\left| a\right| +d}}\sum_{j=0}^{m}(-1)^{2m-j}C_{m}^{j}\Gamma
(N_{j}+\frac{1}{2})2^{N_j+1/2}
\end{equation*}
\begin{equation*}
\times \text{Re}\int_{0}^{\infty }\exp \left( it\frac{\eta ^{\alpha }}{
| z-x| ^\alpha}+\left( \frac{1}{2}N_{j}+\frac{1}{4
}\right) \pi i-\frac{i\pi }{2\alpha }(N_{j}'+1)\right)\end{equation*}
\begin{equation*}
\times
W_{0,N_{j}}(2\eta \exp \left\{ \frac{i\pi (\alpha -1)}{2\alpha }\right\}
)\eta^{N_{j}'}d\eta.
\end{equation*}
Taylor's formula for $\exp ( it\frac{\eta ^{\alpha }}{\left|
z-x\right| ^{\alpha}}) $ yields for the first term, $\forall j\in\leftB 0,m\rightB$,
$$I_\alpha^{j}:=\text{Re}\left\{ \exp \left[ \left( \frac{1}{2}N_{j}+\frac{1}{4}\right) \pi
i-\frac{i\pi }{2\alpha }(N_{j}'+1)\right] \int_{0}^{\infty
}W_{0,N_{j}}(2\eta \exp \left\{ \frac{i\pi (\alpha -1)}{2\alpha }\right\}
)\eta^{N_{j}'}\right. d\eta \bigr\}.$$  
At last, we rotate the contour through the angle $-\frac{\pi (\alpha -1)}{
2\alpha }$. Setting $\xi =\eta \exp \left( \frac{i\pi
(\alpha -1)}{2\alpha }\right) $ we obtain  $I_{\alpha}^{j}=
\text{Re}\left\{ -i(-1)^{j-m}\int_{0}^{\infty }W_{0,N_{j}}(2\xi )\xi
^{N_{j}'}d\xi \right\} =0$.

For $|a| $ and $a_{1}=2m+1$ odd, 
\begin{equation*}
D_{z}^{a }\widetilde{p}^y(t,x,z)=\frac{(-1)^{(| a|
+1)/2}A_{d-2}^{a}}{( 2\pi ) ^{d}| z-x| ^{|
a| +d}} \sum_{j=0}^{m}(-1)^{2m-j}C_{m}^{j}
\end{equation*}
\begin{equation*}
\times\int_{0}^{\infty }dvv^{\left|
a\right| +d-1}\exp \left( -t\frac{v^{\alpha }}{| z-x| ^\alpha}\right) \int_{-1}^{1}(1-\tau ^{2})^{N_{j}-1/2}\tau \sin (v\tau
)d\tau
\end{equation*}
\begin{eqnarray*}
&=&\frac{(-1)^{(|a| +1)/2}A_{d-2}^{a}}{( 2\pi )
^{d}|z-x| ^{|a| +d}}
\sum_{j=0}^{m}(-1)^{2m+1-j}C_{m}^{j}\int_{0}^{\infty }dvv^{| a|
+d-2}\exp \left( -t\frac{v^{\alpha }}{| z-x|^\alpha }\right)\\
&&\times \int_{-1}^{1}(1-\tau ^{2})^{N_{j}-1/2}\tau d(\cos (v\tau ))\\
&=&\frac{(-1)^{(| a| +1)/2}A_{d-2}^{a}}{(2\pi)^{d}| z-x| ^{|a| +d}}
\sum_{j=0}^{m}(-1)^{2m-j}C_{m}^{j}\int_{0}^{\infty }dvv^{| a|
+d-2}\exp \left( -t\frac{v^{\alpha }}{| z-x|^\alpha}\right)\\ 
&&\times \int_{-1}^{1}\cos (v\tau )
\times \left[ (1-2N_{j})\tau^2(1-\tau ^{2})^{N_{j}-3/2}+(1-\tau
^{2})^{N_{j}-1/2}\right] d\tau.
\end{eqnarray*}
The above integrals have the same form as in \eqref{ANC_1A} %\eqref{C_TO_SHOW}
 and can
be estimated similarly.

%\begin{acknowledgements}
%\end{acknowledgements}

\bibliographystyle{alpha}      % basic style, author-year citations
\bibliography{bibli}   % name your BibTeX data base

\end{document}

%% file: macro.tex
\def\I{{\rm{I\!\!I}}}  

\def\cE{{\cal E  }}

\def\F{{\cal F}}

\def\C{{\mathbb{C}}}  
\def\R{{\mathbb{R}} }

\def\N{{\mathbb{N}} }  
\def\E{{\mathbb{E}} }

\def\I{{\mathbb{I}}  }

\def\bint#1^#2{\displaystyle{\int_{#1}^{#2}}}  
\def\bsum#1^#2{\displaystyle{\sum_{#1}^{#2}}}

\def\xdt_#1{X_#1(\Delta t)}

\def\Xtx_#1{X_{#1}^{t,x}}
  
\def\finpreuve{%\begin{flushright}  
$\ \Box$  
}  
\newcommand \A[1]{{\bf (#1)}}

\newtheorem{THM}{Theorem}[section]  
  
\newtheorem{PROP}{Proposition}[section]

\newtheorem{LEMME}{Lemma}[section]  
\newtheorem{REM}{Remark}[section]  
  
\renewcommand{\thesection}{\arabic{section}}  
\newcommand{\mysection}{\setcounter{equation}{0} \section}